\documentclass[12pt]{amsart}
\usepackage{textcase}
\title{$SU(3)$-structures and Special Lagrangian Geometries}
\author{Feng Xu}
\newtheorem{theorem}{Theorem}[section]
\newtheorem{lemma}[theorem]{Lemma}
\newtheorem{proposition}[theorem]{Proposition}
\theoremstyle{definition}
\newtheorem{definition}[theorem]{Definition}

\theoremstyle{remark}
\newtheorem{remark}[theorem]{Remark}

\numberwithin{equation}{section}

\begin{document}

\address{Department of Mathematics, Duke University, Durham, NC 27708}
\email{fxu@math.duke.edu}
\date{}
\begin{abstract}
We generalize Calabi-Yau 3-folds from the special Lagrangian perspective. More precisely, we study $SU(3)$-structures which admit as ``nice" a local special Lagrangian geometry as the flat $\mathbf{C}^3$ or a Calabi-Yau structure does. The underlying almost complex structure may not be integrable. Such $SU(3)$-structures are called {\it admissible}. Among these, we are particularly interested in a subclass of $SU(3)$-structures called {\it nearly Calabi-Yau}. We discuss the local generalities of admissible $SU(3)$-structures and nearly Calabi-Yau structures in Cartan's sense. Examples of nearly Calabi-Yau but non-Calabi-Yau structures as well as other admissible $SU(3)$-structures are constructed from the twistor spaces of self-dual Einstein 4-manifolds. Two classes of examples of complete special Lagrangian submanifolds are found by considering the anti-complex involution of the underlying $SU(3)$-structures. We finally show that the moduli space of compact special Lagrangian submanifolds in a nearly Calabi-Yau manifold behaves in the same way as the moduli in the Calabi-Yau case.
\end{abstract}
\maketitle
\section{Introduction}
The study of special Lagrangian submanifolds was initiated by Harvey and Lawson in their fundamental paper \cite{HarveyLawson}. They proved local existence and constructed many interesting examples in the flat complex spaces. Later on, special Lagrangian geometry received great attention because of its role in mirror symmetry \cite{SYZ}. Since then, interest in special Lagrangian geometry has grown rapidly. 

Three dimensional special Lagrangian geometry is closely related to $SU(3)$-structures. An $SU(3)$-structure on a 6-dimensional manifold $M$ is determined by a real 2-form $\Omega$ and a complex 3-form $\Psi$. A special Lagrangain submanifold can be defined as to anilate $\Omega$ and the imaginary part of $\Psi$. So far, the most thoroughly investigated $SU(3)$-structure is Calabi-Yau which may be defined by the condition $d\Omega=d\Psi=0$. It is well-konwn that the moduli of Calabi-Yau structures on a compact manifold is finite-dimensional. In his interesting program of defining an invariant of Gromov-Witten type using special Lagrangian 3-folds \cite{Joyce}, \cite{Joyce1}, Joyce extended Calabi-Yau to a more general $SU(3)$-structre, the almost Calabi-Yau structure. By doing so, Joyce hoped that one would have enough genericness to simplify the possible singularities of special Lagrangian $3$-folds considerably. The idea of almost Calabi-Yau structures also appeared in the work of Bryant \cite{Bryant4}, and Goldstein \cite{Goldstein}. Almost Calabi-Yau manifolds are K\"ahler. 
 
 In this paper, we study $SU(3)$-structures from the special Lagrangian perspective. We are interested in $SU(3)$-structures which support as many local special Lagrangian submanifolds as the flat $\mathbf{C}^3$ or a Calabi-Yau structure does. This statement can be made precise in the language of exterior differential systems (EDS) \cite{BGC}. Explicitly, we require that the special Lagrangian differential system be involutive. This translates into some first-order conditions on $\Omega$ and $\Psi$. We call the $SU(3)$-structures satisfying these conditions {\it admissible}. 
 
 It is then natural to ask how general admissible $SU(3)$-structures are. An EDS analysis shows that, at least locally, admissible structures are much more general than Calabi-Yau. In particular, the underlying almost complex structures need not be integrable. Although the PDE system defining an admissible $SU(3)$-structure is involutive ``densely", there do exist places where it is not. For example, in some sense, local Calabi-Yau structures are ``singular" points of the moduli of local admissible $SU(3)$-structures. Instead of working with full generality, we consider a smaller class of $SU(3)$-structures called {\it nearly Calabi-Yau}. It is shown that the PDE system for nearly Calabi-Yau is involutive everywhere. This is an indication that the notion of nearly Calabi-Yau may be a better generalizatin of Calabi-Yau. 
  
  To construct complete non-Calabi-Yau admissible examples, we study the twistor spaces of Riemannian 4-manifolds. These spaces carry natural $SU(3)$-structures. If the Riemannian metrics of the base manifolds are self-dual and Einstein, they provide many interesting admissible $SU(3)$-structures including two well-known nearly K\"ahler 6-manifolds. Nearly Calabi-Yau structures appear when the scalar curvatures of the base 4-manifolds are negative. These include the twistor spaces of the hyperbolic 4-space and its various compact quotients.
  
  According to an observation due to Bryant \cite{Bryant2}, the fixed set of an anti-complex involution of an $SU(3)$-structure is special Lagrangian. We construct complete special Lagrangian examples by finding such involutions on the nearly K\"ahler projective 3-space and the twistor space of the hyperbolic 4-space respectively.
  
  In the last section of this paper, we show that the moduli of special Lagrangian submanifolds of a nearly Calabi-Yau manifold behaves as nicely as in the usual Calabi-Yau case.
  
  Finally, it is a great pleasure to thank my advisor, Robert Bryant, for consistent encouragements and helpful discussions, and for his careful revision of an earlier version of this paper. All the errors, if there are any, are however due to myself.

\section{Preliminaries}\label{linearalgebra}
In this section, we set up notational conventions and recall some facts from linear algebra.

\subsection{Metric, orientation and Hodge star}Let $(\mathbf{R}^{2m}, g_0)$ be the standard $2m$ dimensional oriented Euclidean space. Let $\{e_1,e_2,\cdots, e_{2m}\}$ be the standard oriented orthonormal basis. Denote $\{dx_{1},\cdots,dx_{2m}\}$ the dual basis of $(\mathbf{R}^{2m})^*$. There is a Hodge $*$ operator acting on $\Lambda^\bullet (\mathbf{R})^*$ which maps $\Lambda^k$ isomorphically onto $\Lambda^{2m-k}$. In particular, $*$ restricts to an endomorphism on the middle degree piece of $\Lambda^{\bullet}(\mathbf{R})^{2m}$. For example, when $m=2$, \[*: \Lambda^{2}\mathbf{R}^4\rightarrow \Lambda^2\mathbf{R}^4\] has eigenvalues $\pm 1$. The corresponding eigenspaces $\Lambda^2_{\pm}$ are called {\it self-dual} and {\it anti-self-dual} respectively. 
\subsection{Complexification}
For later use we need view $\mathbf{R}^{2m}$ as a subspace of \[\mathbf{C}^{2m}=\mathbf{R}^{2m}\otimes \mathbf{C}.\] We fix a splitting \[\mathbf{C}^{2m}=\mathbf{C}^m\oplus\overline{\mathbf{C}^m}\] by letting $\mathbf{C}^m$ be spanned by \[\{E_{\bar{1}}=\frac{1}{2}(e_1-\sqrt{-1}e_{m+1}),\cdots,E_{\bar{m}}=\frac{1}{2}(e_m-\sqrt{-1}e_{2m})\}\] and $\overline{\mathbf{C}^m}$ spanned by \[\{\overline{E_{\bar{i}}}=\frac{1}{2}(e_i+\sqrt{-1}e_{m+i})\}\] (the reason for the apparent unusual convention of the barred and unbarred indices will be clear later). There is a complex conjugate linear map $\overline{\cdot}$ which maps $\mathbf{C}^{2m}$ onto itself determined by sending $E_{\bar{i}}$ to $\overline{E_{\bar{i}}}$ and vice versa. It interchanges $\mathbf{C}^m$ and $\overline{\mathbf{C}^m}$ and fixes the points in $\mathbf{R}^{2m}$. Dually, \[(\mathbf{R}^{2m})^*\otimes \mathbf{C}=(\mathbf{C}^{2m})^{*}=(\mathbf{C}^{m})^*\oplus \overline{(\mathbf{C}^{m})^*}\] where $(\mathbf{C}^{m})^*$ is complex dual to $\mathbf{C}^m$ and spanned by \[\{dz_i=dx_i+\sqrt{-1}dx_{m+i}\}\] and $\overline{(\mathbf{C}^{m})^*}$ is complex dual to $\overline{\mathbf{C}^m}$ and spanned by \[\{d\overline{z_i}=dx_i-\sqrt{-1}dx_{m+i}\}.\]

We extend the inner product $g_0$ complex bilinearly to $\mathbf{C}^{2m}$, still denoted by $g_0$. It is easy to check that $g_0=dz_1\circ d\overline{z_1}+\cdots+dz_m\circ d\overline{z_m}$. In other words, \[g_0(E_{\bar{i}},E_{\bar{j}})=g_0(\overline{E_{\bar{i}}},\overline{E_{\bar{j}}})=0,\] \[g_0(E_{\bar{i}},\overline{E_{\bar{j}}})=\delta_{\bar{i}j}\] and \[g_0(\overline{E_{\bar{i}}},E_{\bar{j}})=\delta_{i\bar{j}}\]
where $\delta_{\bar{i}i}=1$ and $\delta_{\bar{i}j}=0$ if $i\neq j$ and $\delta_{\bar{i}j}=\delta_{j\bar{i}}$.
\subsection{Linear maps and matrices}
A complex linear map $H:\mathbf{C}^{2m}\rightarrow \mathbf{C}^{2m}$ is identified with a complex matrix 
\begin{equation}\label{matrixrep}\left(\begin{array}{cc}
A_{i\overline{j}}& C_{ij}\\
B_{\overline{ij}} &D_{\overline{i}j}\\

 \end{array}\right)\end{equation}
 by 
 \[H(E_{\bar{i}})=E_{\bar{j}}A_{j\bar{i}}+\overline{E_{\bar{j}}}C_{\overline{ji}}\]
 and 
 \[H(\overline{E_{\bar{i}}})=E_{\bar{j}}B_{ji}+\overline{E_{\bar{j}}}D_{\overline{j}i}.\]
 If $H\in \mathfrak{gl}(2m,\mathbf{R})$ is a linear endomorphism of $\mathbf{R}^{2m}$ we extend it by complex linearity to an endomorphism of $\mathbf{C}^{2m}$, still denoted by $H$. Note that $H(\overline{v})=\overline{H(v)}$ for all $v\in \mathbf{C}^{2m}$. Indeed this charaterizes  $\mathfrak{gl}(2m,\mathbf{R})$. In the matrix representation we have 
 \[\overline{H(E_{\bar{i}})}=\overline{E_{\bar{j}}}\overline{A_{j\bar{i}}}+E_{\bar{j}}\overline{C_{\overline{ji}}}=H(\overline{E_{\bar{i}}}).\]
 Consequently $\overline{A_{j\overline{i}}}=D_{\overline{j}i}$ and $C_{\overline{ji}}=\overline{B_{ji}}$. Thus 
  \[\mathfrak{gl}(2m,\mathbf{R})=\left\{\left(\begin{array}{cc}
                                   A& \overline{B}\\
                                   B&\overline{A}
                                   \end{array}\right)\right\}\subset \mathfrak{gl}(2m,\mathbf{C}).\]
Suppose in addition $H$ preserves $g_0$, i.e., $g_0(Hv,w)+g_0(v,Hw)=0.$ This is equivalent to say that the corresponding matrix satisfies 
 \[\left(\begin{array}{cc}
    A&C\\
    B&D \end{array}\right)
    \left(\begin{array}{cc}
    0 &I\\
    I&0\end{array}\right)+
            \left(\begin{array}{cc}
    0 &I\\
    I&0\end{array}\right)                            
   \left(\begin{array}{cc}
    A&C\\
    B&D \end{array}\right)^t=0.\]
 So the complex orthogonal Lie algebra 
 \[\begin{array}{lll}\mathfrak{so}(2m,\mathbf{C})&=&\left\{\left(\begin{array}{cc}
                                   A& C\\
                                   B&-A^t
                                   \end{array}\right)|B^t+B=C^t+C=0\right\}\\\\
                                   &\subset& \mathfrak{gl}(2m,\mathbf{C}).\end{array}\]
  Hence the real orthogonal Lie algebra is just
  \[\begin{array}{lll}
  \mathfrak{so}(2m,\mathbf{R})&=&\mathfrak{so}(2m,\mathbf{C})\cap\mathfrak{gl}(2m,\mathbf{R})\\\\
  &=&\left\{\left(\begin{array}{cc}
                                   A& \overline{B}\\
                                   B&\overline{A}
                                   \end{array}\right)|B^t+B=0, A^t+\overline{A}=0\right\}.\end{array}\]
It is easy to identify $\mathfrak{u}(m)\subset \mathfrak{so}(2m,\mathbf{R})$ as the set of matrices of the form 
   \[\left(\begin{array}{cc} 
          A&0\\
          0&\overline{A}\end{array}\right)\]where 
          $A^t+\overline{A}=0.$
\subsection{Tensor representation of a linear map}          
The other way to understand $\mathfrak{gl}(2m,\mathbf{C})$ is to view the matrices as tensors of type $(1,1)$. Through the bilinear map $g_0$, $(1,1)$-tensors can be identified with $(2,0)$- or $(0,2)$-tensors. For example if $H$ has its matrix representation (\ref{matrixrep}), then as a $(1,1)$-tensor
  \[H=A_{i\overline{j}}E_{\bar{i}}\otimes dz_{j}+B_{\overline{ij}}\overline{E_{\bar{i}}}\otimes dz_{j}+C_{ij}E_{\bar{i}}\otimes d\overline{z_{j}}+D_{\overline{i}j}\overline{E_{\bar{i}}}\otimes d\overline{z_j},\]
  where the unitary summation convention is understood that summation is implied when an index appears both barred and unbarred.  
  As a $(0,2)$-tensor,
   \[H=\frac{1}{2}(A_{i\overline{j}}d\overline{z_{i}}\otimes dz_{j}+B_{\overline{ij}}dz_i\otimes dz_{j}+C_{ij}d\overline{z_i}\otimes d\overline{z_{j}}+D_{\overline{i}j}dz_i\otimes d\overline{z_j})\]   
while as a $(2,0)$-tensor 
 \[H=2(A_{i\overline{j}}E_{\bar{i}}\otimes \overline{E_{\bar{j}}}+B_{\overline{ij}}\overline{E_{\bar{i}}}\otimes \overline{E_{\bar{j}}}+C_{ij}E_{\bar{i}}\otimes E_{\bar{j}}+D_{\overline{i}j}\overline{E_{\bar{i}}}\otimes E_{\bar{j}}),\]
 because through $g_0$, $E_{\bar{i}}$ corresponds to $\frac{1}{2}d\overline{z_i}$ as a linear functional on $\mathbf{C}^{2m}$, etc.
 If $H\in \mathfrak{so}(2m,\mathbf{R})$, then $D_{\overline{i}j}=-A_{j\overline{i}}$ and $C_{ij}=\overline{B_{\overline{ij}}}=-\overline{B_{\overline{ji}}}.$ Thus the corresponding $(0,2)$-tensor reduces to 
  \[\begin{array}{lll}
  H&=&\frac{1}{2}A_{i\overline{j}}(d\overline{z_{i}}\otimes dz_j-dz_j\otimes d\overline{z_i})\\\\
  &&+\frac{1}{4}B_{\overline{ij}}(dz_i\otimes dz_j-dz_j\otimes dz_i)\\\\
  &&+\frac{1}{4}\overline{B_{\overline{ij}}}(\overline{dz_i\otimes dz_j-dz_j\otimes dz_i}).\end{array}\]
  We adopt the convention that $v\wedge w=v\otimes w-w\otimes v$. Thus $H$ is a two-form
  \[H=\frac{1}{2}A_{i\overline{j}}d\overline{z_i}\wedge dz_j+\frac{1}{4}(B_{\overline{ij}}dz_i\wedge dz_j+\overline{B_{\overline{ij}}}d\overline{z_i}
\wedge d\overline{z_j}).\]
We will use this as the well-known isomorphism $\mathfrak{so}(2m,\mathbf{R})\cong \Lambda^2(\mathbf{R}^{2m})^*.$
Similarly $H$ may be viewed as an element in $\Lambda^2\mathbf{R}^{2m}$, 
\[H=2A_{i\overline{j}}E_{\bar{i}}\wedge \overline{E_{\bar{j}}}+B_{\overline{ij}}\overline{E_{\bar{i}}}\wedge\overline{E_{\bar{j}}}+\overline{B_{\overline{ij}}}E_{\bar{i}}\wedge E_{\bar{j}}.\]
\subsection{Dimension $4$ case}
We now focus our attention on dimension $4$. As mentioned before the Hodge $*$ decomposes 
\[\Lambda^2(\mathbf{C}^4)^*=\Lambda^2_{\mathbf{C}\pm},\] where $\Lambda^2_{\mathbf{C}\pm}=\Lambda^2_{\pm}\otimes\mathbf{C}.$ The set of 2-forms \[\{dz_1\wedge d\overline{z_2}, dz_2\wedge d\overline{z_1},\frac{\sqrt{-1}}{2}(dz_1\wedge d\overline{z_1}-dz_2\wedge d\overline{z_2})\}\] forms a basis for $\Lambda^2_{\mathbf{C}+}$. Similarly \[\{dz_1\wedge dz_2, d\overline{z_1}\wedge d\overline{z_2},\frac{\sqrt{-1}}{2}(dz_1\wedge d\overline{z_1}+dz_2\wedge d\overline{z_2})\}\] is a basis for $\Lambda^2_{\mathbf{C}-}$.

There is a well-known Lie algebra isomorphism between $\mathfrak{so}(4,\mathbf{R})$ and $\mathfrak{su}(2)\oplus \mathfrak{su}(2)$ which we will describe. Let
 \[\begin{array}{lll}\mathfrak{su}(2)_+&=&\left\{\left(\begin{array}{cc}
                                  A& 0\\
                                  0&\overline{A}\end{array}\right)|A^t+\overline{A}=0, tr(A)=0\right\}\\\\
                                  &\subset& \mathfrak{so}(4,\mathbf{R})\end{array}\]
                                  
  and 
  \[\begin{array}{lll}\mathfrak{su}(2)_-&=&\left\{\left(\begin{array}{cc}
                                  \sqrt{-1}aI& \overline{B}\\
                                  B&-\sqrt{-1}aI\end{array}\right)|B^t+B=0, a\in \mathbf{R}\right\}\\\\
                                  &\subset&\mathfrak{so}(4,\mathbf{R}).\end{array}\]   
  It is clear that $\mathfrak{su}(2)_+\cong \mathfrak{su}(2)$ and $\mathfrak{so}(4,\mathbf{R})=\mathfrak{su}(2)_+\oplus\mathfrak{su}(2)_-$. It is also easy to check that $\mathfrak{su}(2)_-$ is a Lie subalgebra and indeed isomorphic to $\mathfrak{su}(2)$. Via the isomorphism $\mathfrak{so}(2m,\mathbf{R})\cong \Lambda^2(\mathbf{R}^{2m})^*$, $\mathfrak{su}(2)_+$ corresponds to $\Lambda^2_+$ and $\mathfrak{su}(2)_-$ corresponds to $\Lambda^2_-$.
  \section{$SU(3)$-structure and special Lagrangian geometry} 
  Let \[g_0=dz_1\circ d\overline{z_1}+\cdots+dz_m\circ d\overline{z_m}\]
      \[\Omega_0=\frac{\sqrt{-1}}{2}(dz_1\wedge d\overline{z_1}+\cdots dz_m\wedge d\overline{z_m})\]
      \[\Psi_0=dz_1\wedge\cdots\wedge dz_m\]    
      be the standard metric, K\"{a}hler form, and holomorphic volume form on $\mathbf{C}^m$. The subgroup of $GL(m,\mathbf{C})$ which preserves these forms is $SU(m)$.
      \subsection{The calibration}
      Harvey and Lawson showed that the $m$-form $\phi_0=
{\rm Re}\Psi_0$ is a calibration on $\mathbf{C}^m$ and that, moreover, an $m$-plane $E\in \mathbf{C}^m$ is $\phi_0$-calibrated if and only if there exists an $A\in SU(m)$ so that $A(E)=\mathbf{R}^m\subset \mathbf{C}^m$. Thus any $\phi_0$ calibrated $E$ satisfies $\Omega_0|_E=0$, i.e., $E$ is an $\Omega_0$-Lagrangian $m$-plane. 
      
      Harvey and Lawson further showed that any $\Omega_0$-Lagrangian $m$-plane $E$ satisfies $\Psi_0|_E=\lambda(E)vol_E$ for some unit complex number $\lambda(E)$. Thus an $\Omega_0$-Lagrangian $m$-plane $E$ is calibrated by $\phi_0$ if and only if $\lambda(E)=1$. For this reason, Harvey and Lawson call the $\phi_0$-calibrated $m$-planes {\it special Lagrangian}. In particular, an $\Omega_0$-Lagrangian $m$-plane $E$ is  $\phi_0$-calibrated with respect to one of its two possible orientations if and only if $\psi_0|_E=0$ where $\psi_0={\rm Im} \Psi_0$. 
      \subsection{SU(m)-structures}
      Let $\pi: P\rightarrow M$ be an $SU(m)$-structure on a manifold of dimension $2m$. The elements of $P_x=\pi^{-1}(x)$ are isomorphisms $u:T_xM\rightarrow \mathbf{C}^m$ and $\pi: P\rightarrow M$ is a principal right $SU(m)$-bundle over $M$ with the right action given by $u\cdot a=a^{-1}\circ u$ for $a\in SU(m)$. Then $P$ defines (and, indeed, is defined by) the metric $g$, the two-form $\Omega$, and the complex $m$-form $\Psi$ defined by 
      \[g_x=u^*(g_0)\]
      \[\Omega_x=u^*(\Omega_0)\]
      \[\Psi_x=u^*(\Psi_0)\]
      for any $u\in P_x$ and $x\in M$.
      Denote the (complex valued) tautological one-forms on $P$ by $\omega=(\omega_1,\omega_2,\cdots,\omega_m)^t$ uniquely determined  by the property
      $\omega_u(X)=u(\pi_*X)$ for any $X\in T_uP$. Then 
      \[\pi^*(\Omega)=\frac{\sqrt{-1}}{2}(\omega_1\wedge\overline{\omega_1}+\cdots+\omega_m\wedge\overline{\omega_m})\]
      and 
      \[\pi^*(\Psi)=\omega_1\wedge\cdots\wedge\omega_m.\]  
   On $P$ there exists a unique $\mathfrak{su}(m)$-valued one-form $\alpha$ called the connection form so that the following structure equations hold
     \begin{equation}\label{strequationSU(m)}
     d\omega_i=-\alpha_{i\overline{j}}\wedge\omega_j+\frac{1}{2}S_{i\overline{jk}}\omega_j\wedge\omega_k+\frac{1}{2}N_{ijk}\overline{\omega_j\wedge\omega_k}+\frac{\sqrt{-1}}{m}(\lambda_k\omega_{k}+\overline{\lambda_k}\omega_k)\wedge\omega_i
     \end{equation}
     where $N_{ijk}=-N_{ikj}$ and $S_{i\overline{jk}}=-S_{i\overline{kj}}$. 
     The quantities $S$, $N$ and $\lambda$ changes tensorially along the fiber and thus are well-defined tensors over $M$. These are called the torsion of the $SU(m)$-structure $P$. For example $N$ is just the famous Nijenhuis tensor and it vanishes if and only if the almost structue implied by the $SU(m)$-structure is integrable. If all the torsions vanish, $P$ is called {\it Calabi-Yau}. The other way to view these torsions is that they are actually components of $d\Omega$ and $d\Psi$ (we omit the details). Thus the $SU(m)$-structure is Calabi-Yau if and only if $d\Omega=d\Psi=0$.
     \subsection{Special Lagrangian submanifolds} A submanifold $L^m\in M^{2m}$ is called {\it special Lagrangian} if it is $\phi={\rm Re}\Psi$ calibrated, i.e., if $L^*(\phi)$ is the volume form. The form $\phi$ is a calibration if it is closed. In this case $L$ is minimal. Assume $L$ is orientable. Then $L$ is special Lagrangian with one of its orientation if and only if $\Omega|_L=\psi|_L=0$ where $\psi={\rm Im}\Psi$. In other words it is an integral manifold of the differential ideal $\mathcal{I}$ generated algebraically by $\Omega$ and $\psi$.
     
     A generic $SU(m)$-structure will not admit any special Lagrangian submanifolds at all, even locally. For example, if \[d\Omega\equiv a\phi\mod (\Omega,\psi)\] for some non-vanishing function $a$, then any special lagrangian submanifold has to anilate $d\Omega$ and hence $\phi$. Thus no special Lagrangian submanifold exists. We are interested to know which $SU(m)$-structures support as many local special Lagrangian submanifolds as the flat $\mathbf{C}^m$ does. If the $SU(m)$-structure is real analytic, so is the ideal $\mathcal{I}$. Now if $d\Omega\in \mathcal{I}$, then we may invoke the Cartan-K\"{a}hler theorem to show that $\mathcal{I}$ is involutive and has the same Cartan characters as the ideal in $\mathbf{C}^m$. Here we need not care about whether or not $d({\rm Im}\Psi)$ is in $\mathcal{I}$ because $d({\rm Im}\Psi)$ is of degree $m+1$ and hence vanishes automatically on any $m$-dimensional submanifold. To analyze further we need distinguish when $m>3$ and when $m=3$. 
     
     When $m>3$ we are requiring the existence of a 1-form $\theta$ such that $d\Omega=\theta\wedge\Omega$. This translates into some condition on the torsion. We will leave it for the reader.
     We are mainly interested in $m=3$, where we have more flexibility.
    \subsection{Admissible $SU(3)$-structures} We make the following definition.
\begin{definition}\label{Admissible $SU(3)$-structure}(Admissible $SU(3)$-structures)
 An $SU(3)$-structure $(\Omega, \Psi)$ on $M$ is called {\it admissible} if there exist a $1$-form $\theta$ and a real function $a$ such that
 \begin{equation}\label{admissiblecondition}
     d\Omega=\theta\wedge\Omega+a\psi.\end{equation}
\end{definition}
     To unravel this condition, let us first rewrite the structure equations (\ref{strequationSU(m)}) adapted to $m=3$. Let \[N_{i\overline{j}}=\frac{1}{2}\overline{\epsilon_{jkl}}N_{ikl}\] and \[S_{ij}=\frac{1}{2}\epsilon_{jkl}S_{i\overline{kl}}\] 
    where $\epsilon_{ijk}$ is the signature of the permutation $(ijk)$. Note these are invertible and the inverse transformations are 
     $N_{ijk}=\epsilon_{jkl}N_{i\overline{l}}$ and 
     $S_{i\overline{jk}}=\overline{\epsilon_{jkl}}S_{il}$. 
     We rewrite the structure equations as
    \begin{equation}\label{strequationSU(3)}     d\omega_i=-\alpha_{i\overline{j}}\wedge\omega_j+\frac{1}{2}S_{ij}\overline{\epsilon_{jkl}}\omega_k\wedge\omega_l+\frac{1}{2}N_{i\overline{j}}\epsilon_{jkl}\overline{\omega_k\wedge\omega_l}+\frac{\sqrt{-1}}{3}(\lambda_k\overline{\omega_k}+\overline{\lambda_k}\omega_k)\wedge\omega_i.
         \end{equation}
     
  Using this we can easily compute that 
  \[\pi^*(d\Omega)=\frac{\sqrt{-1}}{2}(N_{i\overline{i}}\overline{\Psi}-\overline{N_{i\overline{i}}}\Psi)+\frac{\sqrt{-1}}{4}(S_{il}\overline{\epsilon_{ljk}}\overline{\omega_i}\wedge\omega_j\wedge\omega_k-\overline{S_{il}}\epsilon_{ljk}\omega_i\wedge\overline{\omega_j\wedge\omega_k}).      \]   
    Suppose condition (\ref{admissiblecondition}) is satisfied with 
    \[\theta=u_{\overline{i}}\omega_i+\overline{u_{\overline{i}}}\overline{\omega_i}.\]
    Then it also holds that
    \[d\Omega=\frac{\sqrt{-1}}{4}(u_{\overline{i}}\delta_{\overline{j}k}-u_{\overline{j}}\delta_{\overline{i}k})\omega_i\wedge\omega_j\wedge\overline{\omega_k}-\frac{\sqrt{-1}}{4}(\overline{u_{\overline{i}}}\delta_{\overline{k}j}-\overline{u_{\overline{j}}}\delta_{\overline{k}i})\overline{\omega_i\wedge\omega_j}\wedge\omega_k+a\psi.\] 
    Comparing the two expressions for $d\Omega$ we get
    \[\overline{\epsilon_{ljk}}S_{il}=u_{\overline{j}}\delta_{i\overline{k}}-u_{\overline{k}}\delta_{i\overline{j}}\]
    and
    $a=N_{i\overline{i}}$. 
  We summarize the discussion so far as follows
  \begin{lemma}\label{involutivity}
  The ideal $\mathcal{I}$ is involutive and every $\Omega$-Lagrangian analytic $2$-submanifold in $M^6$ can be thickened uniquely to a special Lagrangian submanifold if there exist a (necessarily unique) connection 
  $\alpha$ so that 
  \begin{equation}\label{admissiblestrequation}
d\omega_i=-\alpha_{i\overline{j}}\wedge\omega_j+\beta\wedge\omega_i+\frac{1}{2}N_{i\overline{j}}\epsilon_{jkl}\overline{\omega_k\wedge\omega_l}\end{equation}
 where $\beta$ is a complex $1$-form and the trace of the Nijenhuis tensor, $N_{i\overline{i}}$, is real.    
  \end{lemma}
\begin{remark}
The same result for $\mathbf{C}^3$ was shown by Harvey and Lawson in \cite{HarveyLawson}. Lemma \ref{involutivity} says an admissible $SU(3)$-manifold supports as nice a local special Lagrangian geometry as $\mathbf{C}^3$ does. This is the best situation one can hope. 
 \end{remark}
\begin{remark}
  It is easy to see that then condition (\ref{admissiblestrequation}) in Lemma \ref{involutivity} essentially gives an equivalent definition of admissible $SU(3)$-structures. For this reason, we will also call an $SU(3)$-structure satisfying (\ref{admissiblestrequation}) admissible.
 \end{remark} 
  If, in addition, we require $d\psi=0$, more interesting conditions appear. Using the structure equations we compute
   \[d\psi=\frac{\sqrt{-1}}{8}\overline{\epsilon_{ijk}}\epsilon_{lpq}(N_{i\overline{l}}-\overline{N_{l\overline{i}}})\overline{\omega_{p}\wedge\omega_q}\wedge\omega_j\wedge\omega_k+\frac{\sqrt{-1}}{2}(\overline{\lambda_{\overline{l}}}\overline{\omega_{l}}\wedge\Psi+\lambda_{\overline{l}}\omega_l\wedge\overline{\Psi}).\]
  Thus $\lambda=0$ and $N_{i\overline{j}}-\overline{N_{j\overline{i}}}=0$. Out of these we pick a class of special interest and call it {\it nearly Calabi-Yau}. 
  \begin{definition}\label{nearlycalabiyau} (Nearly Calabi-Yau)
  An $SU(3)$-structure $(M^6,\Omega, \Psi)$ is called {\it nearly Calabi-Yau} if 
  \[d\Omega=d({\rm Im}\Psi)=0.\]  
  \end{definition} 
  Thus, the difference from Calabi-Yau is we do not require $d({\rm Re}\Psi)=0$. We will call it {\it strictly nearly Calabi-Yau} if $d({\rm Re}\Psi)\neq 0$.
  In terms of structure equations we have
  \begin{proposition}[Alternative Definition of Nearly Calabi-Yau]
  An $SU(3)$-structure is nearly Calabi-Yau if and only if there exists a (necessarily unique) connection $\alpha$ so that \begin{equation}\label{nearlycalabiyaustreqn} d\omega_i=-\alpha_{i\overline{j}}\wedge\omega_j+\frac{1}{2}N_{i\overline{j}}\epsilon_{jkl}\overline{\omega_k\wedge\omega_l}, \end{equation}
where $N_{i\overline{j}}-\overline{N_{j\overline{i}}}=0$ and $tr(N)=N_{i\overline{i}}=0$. 
  \end{proposition}
 In other words the Nijenhuis tensor is Hermitian symmetric and trace free and no other torsion exists. 
 \subsection{Examples}
\subsubsection{Calabi-Yau}
This is the case of most interest so far. In the next section we will give strictly nearly Calabi-Yau examples by studying twistor spaces of certain Riemannian $4$-manifolds.
\subsubsection{Nearly K\"{a}hler} 
The fundamental forms $\Omega$ and $\Psi$ satisfy 
\[d\Omega=3c{\rm Im}\Psi\] and \[d\Psi=2c\Omega^2.\]
It is well-known that the underlying metric is Einstein by \cite{De-Ka}. It also follows that such structures are real analytic, in, say coordinates harmonic for the metric. The ideal $\mathcal{I}$ is clearly differentially closed. Hence the almost special Lagrangian geometry of nearly K\"{a}hler $6$ manifolds is well-behaved locally. 

When $c=0$, the structure is actually Calabi-Yau. When $c\neq 0$, we can scale the forms to reduce the structure to the case $c=1$. The underlying almost complex structure structure is non-integrable. In fact, the Nijenhuis tensor $N_{i\overline{j}}$ is the identity matrix. Nearly K\"{a}hler but non-Calabi-Yau Examples include $S^6$ with the standard metric and almost complex structure, $S^3\times S^3$, the flag manifold $SU(3)/T^2$ and the projective space $\mathbf{CP}^3$ (with an unusual almost complex structure, however).
\begin{remark}
It should be cautioned that a nearly Calabi-Yau manifold, unless it is Calabi-Yau, is NOT nearly K\"{a}hler.
\end{remark}
\subsection{Generalities}
Calabi-Yau and nearly K\"{a}hler provide first examples of admissible $SU(3)$-structures. However, we are about to show that they are only a ``closed" subset of the moduli of local admissible $SU(3)$-structures. For this we need study the local generality of admissible $SU(3)$-structures.  

Let $p:\mathcal{F}\rightarrow  M^6$ be the total coframe bundle of $M^6$. Thus $\mathcal{F}$ is a principal $GL(6,\mathbf{R})$-bundle over $M^6$. The fiber $\mathcal{F}_x$ over $x$ consists of the linear isomorphisms $u:T_xM\rightarrow \mathbf{R}^6$. We consider the quotient $\mathcal{F}/SU(3)$  which is $34(=6+6^2-8)$ dimensional and projects onto $M$ with fibers diffeomorphic to $GL(6,\mathbf{R})/SU(3)$. An $SU(3)$-structure over $M$ may be regarded as a section of $\overline{p}:\mathcal{F}/SU(3)\rightarrow M^6$. In fact, if $P$ is an $SU(3)$-structure, then $P$ is a subbundle of $\mathcal{F}$ and every fiber $P_x$ determines a unique $SU(3)$ orbit of $\mathcal{F}_x$ and thus a section of $\mathcal{F}/SU(3)\rightarrow M^6$. Conversely, if $\sigma: M\rightarrow \mathcal{F}/SU(3)$ is such a section, let $P$ be the preimage of $\sigma(M)$ under the projection $\mathcal{F}\rightarrow \mathcal{F}/SU(3)$. Then $P$ is the needed $SU(3)$-structure.  
\begin{remark}
There is a more concrete realization of $\mathcal{F}/SU(3)$. Let $\Lambda^2_+M\oplus\Lambda^3_+M$ be the subbundle of  $\Lambda^2 M \oplus\Lambda^3 M$ consiting of the pairs of {\it positive forms} $(\rho^2,\rho^3)$ in the sense that there exists a linear isomorphism $u:T_xM\rightarrow \mathbf{R}^6$ such that $\rho^2=u^*\Omega_0$ and $\rho^3=u^*{\rm Re}(\Psi_0)$. We clearly have a projection $\mathcal{F}\rightarrow \Lambda^2_+M\oplus\Lambda^3_+M$. Since the isotropy group of $(\Omega_0,{\rm Re}\Psi_0)$ is $SU(3)\subset GL(6,\mathbf{R})$, this is indeed a principal $SU(3)$-bundle. For similar discussions concerning $G_2$-structures and $Spin(7)$-structures, consult Bryant's work \cite{Bryant1}. In fact, the work directly inspired the discussion in this section. 
\end{remark}
We will write structures equations for $\mathcal{F}$. For notational conventions on linear algebra see Section \ref{linearalgebra}. 
Let $(\omega_1,\omega_2,\omega_3,\overline{\omega_1},\overline{\omega_2},\overline{\omega_3})$ denote complexified tautological one-forms. Thus, for example, $\omega_{1u}=u^*(dz_1)\circ p_*.$
We fix a $\mathfrak{gl}(6,\mathbf{R})$-valued connection form 
 \[\left(\begin{array}{cc}
 \alpha+\frac{\sqrt{-1}}{3}\mu+\kappa &\beta\\
 \overline{\beta}& \overline{\alpha}-\frac{\sqrt{-1}}{3}\mu+\overline{\kappa}\\
 \end{array}\right)\]on $\mathcal{F}$ where $\alpha$ takes value in $\mathfrak{su}(3)$, $\beta$ is $\mathfrak{gl}(3,\mathbf{C})$-valued, $\mu$ is a real form, and $\kappa$ satisfies
 \[\overline{\kappa}^t=\kappa.\] 
 We have the following structure equations on $\mathcal{F}$:
 \begin{equation}
 d\left(\begin{array}{c}
    \omega\\
    \overline{\omega}\end{array}\right)=
    -\left(\begin{array}{cc}
 \alpha+\frac{\sqrt{-1}}{3}\mu+\kappa &\beta\\
 \overline{\beta}& \overline{\alpha}-\frac{\sqrt{-1}}{3}\mu+\overline{\kappa}\\
 \end{array}\right)\wedge
 \left(\begin{array}{c}
    \omega\\
    \overline{\omega}\end{array}\right).
  \end{equation}
 Note that the forms $\Omega=\frac{\sqrt{-1}}{2}(\omega_1\wedge\overline{\omega_1}+\omega_2\wedge\overline{\omega_2}+\omega_3\wedge\overline{\omega_3})$ and $\Psi=\omega\wedge\omega_2\wedge\omega_3$ are $SU(3)$ invariant on $\mathcal{F}$, so they descend to $\mathcal{F}/SU(3)$. We use the same letters to denote the forms on $\mathcal{F}/SU(3)$ and let $\psi={\rm Im}\Psi$ and $\phi={\rm Re}\Psi$. We will use Cartan-K\"{a}hler machinary to study local admissible $SU(3)$-structures and nearly Calabi-Yau structures. For background material on exterior differential systems, see the standard text \cite{BGC}.
 \subsubsection{Generalities of admissible $SU(3)$-structures} We introduce a new manifold $\mathcal{M}=(\mathcal{F}\times \mathbf{C}^3)/SU(3)\times \mathbf{R}$, where $SU(3)$ acts on $\mathbf{C}^3$ in the obvious way. We use $(u_{\overline{1}},u_{\bar{2}},u_{\bar{3}})$ as the coordinate on $\mathbf{C}^3$ and $a$ as the coordinate on $\mathbf{R}$. $\mathcal{M}$ is a vector bundle of rank $7$ over $\mathcal{F}/SU(3)$. Let $\theta=u_{\bar{i}}\omega_i+\overline{u_{\bar{i}}}\overline{\omega_i}$. Then $\theta$ is another well-defined differential form on $\mathcal{M}$ besides $\Omega$ and $\Psi$. On $\mathcal{M}$ define a differential ideal
 \[
 \begin{array}{rrl}
 \mathbf{I}&=&\langle\Pi_3=d\Omega-\theta\wedge\Omega-a\psi\rangle_{diff}\\
           &=&\langle\Pi_3=d\Omega-\theta\wedge\Omega-a\psi,       \Pi_4=d\theta\wedge\Omega+(da-a\theta)\wedge\psi+ad\psi\rangle_{alg}.\end{array}\]
  We are interested in $6$-dimensional (local) integral manifolds of this ideal which are also  local sections of $\mathcal{M}\rightarrow M^6$. Such a section pulls back $\Omega$ and $\Psi$ to $M$ which satisfies the condition (\ref{admissiblecondition}) and thus defines an admissible $SU(3)$-structure. A section satisfies $\sqrt{-1}\Psi\wedge\overline{\Psi}\neq 0$. Conversely a $6$-dimensional submanifold of $\mathcal{M}$ on which $\sqrt{-1}\Psi\wedge\overline{\Psi}\neq 0$ is locally a section. Hence we will consider the integral manifolds of $\mathbf{I}$ with the {\it independence condition} $\sqrt{-1}\Psi\wedge\overline{\Psi}\neq 0.$ 
  \begin{theorem}
  The differential system $(\mathbf{I},\sqrt{-1}\Psi\wedge\overline{\Psi}\neq 0)$ on the dense open set $\mathcal{M}\setminus \{a= 0\}$ is involutive with Cartan characters $(s_0,s_1,s_2,s_3,s_4,s_5,s_6)=(0,0,1,3,6,10,15).$
  \end{theorem}
  \begin{proof}
  Since the system contains no forms of degree $2$ or less, we have $c_0=c_1=0$. On the other hand, $c_1=1$ and $c_6=35$. We need compute the remaining three characters $c_3$, $c_4$ and $c_5$. For this we pass up to $\mathcal{F}\times \mathbf{C}^3\times \mathbf{R}$. For effective computations we set
 \[Da=da-a\theta,\]
 and  \[Du_{\bar{i}}=du_{\bar{i}}-u_{\bar{j}}\alpha_{j\bar{i}}-u_{\bar{j}}\kappa_{j\bar{i}}-\frac{\sqrt{-1}}{3}u_{\bar{i}}\mu-\overline{u_{\bar{j}}}\overline{\beta_{ji}}.\]
   Relative to the projection $\mathcal{F}\times \mathbf{C}^3\times \mathbf{R}\rightarrow \mathcal{M}$, the forms $\omega, \kappa, \beta, \mu, Da, Du$ form a basis for semibasic $1$-forms. In terms of these forms we have
    \[\begin{array}{ccl}
    \Pi_3&=&-\frac{\sqrt{-1}}{2}\kappa_{i\bar{j}}\wedge\omega_j\wedge\overline{\omega_i}+\frac{\sqrt{-1}}{2}\overline{\kappa_{i\bar{j}}}\wedge\overline{\omega_j}\wedge\omega_i\\\\
    &&-\frac{\sqrt{-1}}{4}(\beta_{ij}-\beta_{ji})\wedge\overline{\omega_j}\wedge\overline{\omega_i}
    +\frac{\sqrt{-1}}{4}(\overline{\beta_{ij}}-\overline{\beta_{ji}})\wedge\omega_j\wedge\omega_i\\\\
    &&-\frac{\sqrt{-1}}{2}(u_{\bar{i}}\omega_i+\overline{u_{\bar{i}}}\overline{\omega_i})\wedge\omega_j\wedge\overline{\omega_j}\\\\
   && -\frac{\sqrt{-1}}{2}a(\overline{\omega_1\wedge\omega_2\wedge\omega_3}-\omega_1\wedge\omega_2\wedge\omega_3)
    \end{array}\]  
 and 
 \[\begin{array}{ccl}
 \Pi_4&=&\frac{\sqrt{-1}}{2}(Du_{\bar{i}}\wedge\omega_i+\overline{Du_{\bar{i}}}\wedge\overline{\omega_i})\wedge\omega_j\wedge\overline{\omega_j}\\\\
     &&+\frac{\sqrt{-1}}{2}(Da-a\kappa_{i\bar{i}})\wedge(\overline{\omega_1\wedge\omega_2\wedge\omega_3}-\omega_1\wedge\omega_2\wedge\omega_3)\\\\
  &&-\frac{a}{2}\mu\wedge(\overline{\omega_1\wedge\omega_2\wedge\omega_3}+\omega_1\wedge\omega_2\wedge\omega_3)\\\\
 &&-\frac{\sqrt{-1}}{4}a\epsilon_{ijk}\overline{\beta_{il}}\wedge\omega_l\wedge\overline{\omega_j\wedge\omega_k}\\\\
 &&+\frac{\sqrt{-1}}{4}a\overline{\epsilon_{ijk}}\beta_{il}\wedge\overline{\omega_l}\wedge\omega_j\wedge\omega_k.
 \end{array}
 \]   
 A six dimensional subspace $E_6$ of the tangent plane on which $\Psi\wedge\overline{\Psi}\neq 0$ is defined by the following relations
 \begin{equation}\label{parameters}\left\{\begin{array}{ccl}
   \kappa_{i\overline{j}}&=&A_{i\overline{j}k}\overline{\omega_k}+\overline{A_{j\bar{i}k}}\omega_k;\\\\
 \beta_{ij}&=&B_{ijk}\overline{\omega_k}+C_{ij\overline{k}}\omega_k,\\\\
 Du_{\bar{i}}&=&U_{\overline{ij}}\omega_{j}+U_{\overline{i}j}\overline{\omega_j},\\\\
 Da&=&a_{\bar{i}}\omega_i+\overline{a_{\bar{i}}}\overline{\omega_i},\\\\
 \mu&=&b_{\bar{i}}\omega_i+\overline{b_{\bar{i}}}\overline{\omega_i},   
 \end{array}\right.   \end{equation}
 where $A_{i\overline{j}k}, B_{ijk}, C_{ij\overline{j}}, U_{\overline{ij}}, U_{\overline{i}j}$ and $a_{\bar{i}}, b_{\bar{i}}$ are free parameters. In order that $E_6$ be an integral element of $\mathbf{I}$, it must anilate $\Pi_3$ and $\Pi_4$. This amounts to the following equations on the parameters in (\ref{parameters}),
 \begin{equation}\label{equationE6}
 \left\{\begin{array}{rcl}
 \epsilon_{ijk}\overline{B_{ijk}}&=&a,\\\\
 2(A_{k\bar{j}i}-A_{i\bar{j}k})+C_{ik\bar{j}}-C_{ki\bar{j}}+\overline{u_{\bar{i}}}\delta_{\bar{j}k}-\overline{u_{\bar{k}}}\delta_{\bar{j}i}&=&0,\\\\
 \sqrt{-1}(U_{\overline{ij}}\epsilon_{jik}+\overline{a_{\bar{k}}}-C_{ik\bar{i}})-a(\overline{b_{\bar{k}}}-\sqrt{-1}A_{i\overline{i}k})&=&0,\\\\
 -(U_{\overline{i}j}\delta_{k\overline{l}}+U_{\overline{l}k}\delta_{j\overline{i}}-U_{\overline{l}j}\delta_{k\overline{i}}-U_{\overline{i}k}\delta_{j\overline{l}})&&\\\\
 +(\overline{U_{\bar{k}l}}\delta_{j\bar{i}}-\overline{U_{\bar{j}l}}\delta_{k\bar{i}}- \overline{U_{\bar{k}i}}\delta_{j\bar{l}}+ \overline{U_{\bar{j}i}}\delta_{k\bar{l}})& &
 \\\\-a\overline{\epsilon_{pil}}B_{pkj}+a\overline{\epsilon_{pil}}B_{pjk}+a\epsilon_{pjk}\overline{B_{pli}}-a\epsilon_{pjk}\overline{B_{pil}}&=&0.\\\\
 \end{array}\right.
 \end{equation}
By inspection, we have $35=(2+3\times 3\times 2+3\times 2+3\times 3)$ linearly independent affine equations in (\ref{equationE6}) (note that the last equations are real, while the others are complex) . The solution space is smooth, even where $a=0$. We pick $E_5=\text{span}\{e_1,e_2,e_3,e_4,e_5\}\subset E_6$ where $e_1$ is dual to ${\rm Re}(\omega_1)$ and $e_4$ is dual to ${\rm Im}(\omega_1)$, etc. First note that 
\[c_5\leq\left(\begin{array}{c}
       5\\
       2 \end{array}\right)+\left(\begin{array}{c}
       5\\
       3 \end{array}\right)=20.           \]
 We will show that the equality holds, i.e., the polar equations of $E_5$ has the largest possible rank $20$. It then follows that if we pick any flag $E_3\subset E_4\subset E_5$ we have \[c_3=\left(\begin{array}{c}
       3\\
       2 \end{array}\right)+\left(\begin{array}{c}
       3\\
       3 \end{array}\right)=4\]
and
   \[c_4=\left(\begin{array}{c}
       4\\
       2 \end{array}\right)+\left(\begin{array}{c}
       4\\
       3 \end{array}\right)=10. \]
  Since $c_0+c_1+c_2+c_3+c_4+c_5=1+4+10+20=35$ we apply Cartan's test to finish the proof.
 
  The verification that the polar equations of $E_5$ have rank $20$ is a lengthy linear algebra exercise. First, by translating $\kappa,\beta, Du, Da, \mu$ we may assume these forms vanish on $E_6$ since $\Pi_3$ and $\Pi_4$ are affine linear in these forms. Now the rank of polar equations are the number of linearly independent forms in 
  $\{(e_{i}\wedge e_j)\lrcorner\Pi_3, (e_i\wedge e_j\wedge e_k)\lrcorner\Pi_4\}$. We omit the messy details but only point out the following facts (an unsatisfied reader may consult the computations in the proof of Theorem (\ref{nearlycalabiyaugenerality}) and make necessary modifications by himself). The forms $\{(e_{i}\wedge e_j)\lrcorner\Pi_3\}$ pick out $10$ linearly independent forms from linear combinations of ${\rm Re}(\kappa_{i\overline{j}})$, ${\rm Im} (\kappa_{i\bar{j}})$, ${\rm Re} (\beta_{ij}-\beta_{ji})$ and ${\rm Im}(\beta_{ij}-\beta_{ji})$. The $6$ forms $(e_i\wedge e_{i+3}\wedge e_j)\lrcorner \Pi_4$ pick out real and imaginary parts of $\{Du_{\bar{j}}+\text{linear combinations of }\beta\}_{j=1}^{2}\bigcup{\{{\rm Re}(Du_{\bar{3}}-a\beta_{21}), {\rm Re}(Du_{\bar{3}}+a\beta_{12})\}}$ . The $4$ forms $(e_{i}\wedge e_j\wedge e_k)\lrcorner \Pi_4$ with $i\in \{1,4\}$, $j\in\{2,5\}$ and $k=3$ give us non-degenerate linear combinations of the forms $Da-a\sum\kappa_{i\overline{i}}$, $a\mu$, $a{\rm Re}(\beta_{ii})$ and $a{\rm Im}(\beta_{ii})$.  These three classes of equations are clearly independent from each other if we assume $a\neq 0$.    
  \end{proof}
 \begin{remark} There does not exist any regular flag over the locus $\{a=0\}$ for simple reasons. When $a=0$, only $7$ independent forms $Du_{\bar{i}}$ and $Da$ in $\Pi_4$ could contribute to the polar equations. Thus $c_5\leq 17$.
 \end{remark} 
  \begin{remark}The last nonzero character is $s_6=15$. Modulo diffeomorphisms, which depend on $6$ functions of $6$ variables, we still have $9$ functions of $6$ variables of local gererality of admissible $SU(3)$-structures. Both local Calabi-Yau and nearly K\"{a}hler structures depend on $2$ functions of $5$ variables. Thus local admissible $SU(3)$-structures are much more general than Calabi-Yau and nearly K\"{a}hler.  
  \end{remark}
\subsubsection{Generalities of nearly Calabi-Yau} Nearly Calabi-Yau is a subclass of admissible $SU(3)$-structures. At the first thought, one would expect nearly Calabi-Yau is less general than an admissible $SU(3)$-structure. We will show this is indeed the case. 
Now the differential system $\mathbb{I}$ is defined on $\mathcal{F}/SU(3)$ and generated algebraically by the $3$-form $d\Omega$ and the $4$-form $d\psi$ with the independence condition $\sqrt{-1}\Psi\wedge\overline{\Psi}$. This system is better-behaved than $\mathbf{I}$ for admissible $SU(3)$-structures in that it is involutive on the whole $\mathcal{F}/SU(3)$.  
\begin{theorem}\label{nearlycalabiyaugenerality}
The differential system $\mathbb{I}$ on $\mathcal{F}/SU(3)$ is involutive with Cartan characters $(s_0,s_1,s_2,s_3,s_4,s_5,s_6)=(0,0,1,3,6,9,9).$
 \end{theorem} 
 \begin{proof}Since the system contains no forms of degree $2$ or less, we have $c_0=c_1=0$. Moreover, it is easy to see $c_2=1$. To use Cartan's test, we need compute the other $3$ characters $c_3,c_4$ and $c_5$ and the codimension of the space of $6$-dimensional integral elements. For this we pass up to $\mathcal{F}$ where 
 \[\begin{array}{ccl}
 d\Omega&=&-\frac{\sqrt{-1}}{2}\kappa_{i\bar{j}}\omega_j\wedge\overline{\omega_i}+\frac{\sqrt{-1}}{2}\overline{\kappa_{i\bar{j}}}\overline{\omega_j}\wedge\omega_i\\\\
       &&-\frac{\sqrt{-1}}{2}\beta_{ij}\wedge\overline{\omega_j}\wedge\overline{\omega_i}+\frac{\sqrt{-1}}{2}\overline{\beta_{ij}}\wedge\omega_j\wedge\omega_i,
       \end{array} \]
 and
 \[\begin{array}{ccl}
 d\psi&=&-\frac{1}{2}(\mu+\sqrt{-1}\kappa_{i\overline{i}})\wedge\overline{\omega_1\wedge\omega_2\wedge\omega_3}- \frac{1}{2}(\mu-\sqrt{-1}\kappa_{i\overline{i}})\wedge\omega_1\wedge\omega_2\wedge\omega_3\\\\   
    &&-\frac{\sqrt{-1}}{4}\epsilon_{ijk}\overline{\beta_{il}}\wedge\omega_l\wedge\overline{\omega_j}\wedge\overline{\omega_k}+\frac{\sqrt{-1}}{4}\overline{\epsilon_{ijk}}\beta_{il}\wedge\overline{\omega_l}\wedge\omega_j\wedge\omega_k.
    \end{array}\]    
 A $6$-dimensional integral element $E_6$ on which $\sqrt{-1}\Psi\wedge\overline{\Psi}\neq 0$
is parametrized by the equations in (\ref{parameters}) for $\kappa, \beta$ and $\mu$ but now the quantities $A, B$ and $b$ satisfy the following equations 
\begin{equation}\left\{
\begin{array}{rcl}
\overline{\epsilon_{ijk}}B_{ijk}&=&0\\\\
2A_{k\overline{j}i}-2A_{i\overline{j}k}+C_{ik\overline{j}}-C_{ki\overline{j}}&=&0\\\\
C_{ik\overline{i}}-A_{i\overline{i}k}-\sqrt{-1}\overline{b_{\bar{k}}}&=&0\\\\
-\overline{\epsilon_{pil}}B_{pkj}+\overline{\epsilon_{pil}}B_{pjk}+\epsilon_{pjk}\overline{B_{pli}}-\epsilon_{pjk}\overline{B_{pil}}&=&0.
\end{array}\right.
\end{equation}
 The last equations are real while the others are all complex. Moreover, the last equations imply the imaginary part of the first equation. Thus the total rank of these linear equations is $1+3\times 3\times 2+3\times 2+3\times 3=34$. The forms $\rm{Re}(\omega_i)$ and ${\rm Im} (\omega_i)$ restrict to $E_6$ to be a dual basis. Let $\{e_1, e_2,e_3,e_4,e_5,e_6\}$ be the basis of $E_6$ for which $e_1$ is dual to $\rm{Re}(\omega_1)$ and $e_4$ is dual to ${\rm Im}(\omega_1)$, etc. Again by translating we may assume $A=B=b=0.$ Let $E_3=\text{span}\{e_1,e_2,e_3\}$, $E_4=\text{span}\{e_1,e_2,e_3,e_4\}$, and $E_5=\text{span}\{e_1,e_2,e_3,e_4,e_5\}$.
 
 The polar space of $E_3$ consists of vectors anilating the following $1$-forms
 \begin{equation}\label{polarE3}\left\{
 \begin{array}{rcl}
 (e_1\wedge e_2)\lrcorner d\Omega&=&-2{\rm Im}(\kappa_{1\bar{2}})-{\rm Im}(\beta_{12}-\beta_{21})\\
 (e_1\wedge e_3)\lrcorner d\Omega&=&-2{\rm Im}(\kappa_{1\bar{3}})-{\rm Im}(\beta_{13}-\beta_{31})\\
 (e_2\wedge e_3)\lrcorner d\Omega&=&-2{\rm Im}(\kappa_{2\bar{3}})-{\rm Im}(\beta_{23}-\beta_{32})\\
 (e_1\wedge e_2\wedge e_3)\lrcorner d\psi&=& \mu+\rm{Im}(\beta_{11}+\beta_{22}+\beta_{33})\\
  \end{array}\right\}.
 \end{equation} 
 Consequently $c_3=4$ and $s_3=3$.
 
 The polar space of $E_4$ consists of vectors anilating the forms in (\ref{polarE3}) as well as the following forms
 \begin{equation}\label{polarE4}\left\{
 \begin{array}{rcl}
 (e_1\wedge e_4)\lrcorner d\Omega&=&-2\kappa_{1\bar{1}}\\
 (e_2\wedge e_4)\lrcorner d\Omega&=&-2\rm{Re}(\kappa_{1\bar{2}})-\rm{Re}(\beta_{12}-\beta_{21})\\
 (e_3\wedge e_4)\lrcorner d\Omega&=&-2\rm{Re}(\kappa_{1\bar{3}})-\rm{Re}(\beta_{13}-\beta_{31})\\
 (e_1\wedge e_2\wedge e_4)\lrcorner d\psi&=&-2{\rm Re}(\beta_{31})\\
 (e_1\wedge e_3\wedge e_4)\lrcorner d\psi&=&2{\rm Re}(\beta_{21})\\
 (e_2\wedge e_3\wedge e_4)\lrcorner d\psi&=&-\sum_{i}\kappa_{i\overline{i}}+{\rm Re}(-\beta_{11}+\beta_{22}+\beta_{33}) 
 \end{array}
 \right\}.
 \end{equation}
 These forms are independent among themselves and also independent from forms in (\ref{polarE3}). Thus $s_4=6$.
 The polar space for $E_5$ consists of vector anilating forms in (\ref{polarE3}), (\ref{polarE4}) as well as the following forms 
 \begin{equation}\label{polarE5}\left\{
 \begin{array}{rcl}
 (e_1\wedge e_5)\lrcorner d\Omega &=&-2{\rm Re}(\kappa_{1\bar{2}})+\rm{Re}(\beta_{12}-\beta_{21})\\
 (e_2\wedge e_5)\lrcorner d\Omega &=&-2\kappa_{2\bar{2}}\\
 (e_3\wedge e_5)\lrcorner d\Omega &=&-2{\rm Re}(\kappa_{3\bar{2}})+{\rm Re}(\beta_{32}-\beta_{23})\\
 (e_4\wedge e_5)\lrcorner d\Omega &=&-2{\rm Im}(\kappa_{1\bar{2}})+{\rm Im}(\beta_{12}-\beta_{21})\\
 (e_1\wedge e_2\wedge e_5)\lrcorner d\psi&=&-2{\rm Re}(\beta_{32})\\
 (e_1\wedge e_3\wedge e_5)\lrcorner d\psi&=&-\sum_{i}\kappa_{i\overline{i}}+{\rm Re}(-\beta_{11}+\beta_{22}-\beta_{33})\\
 (e_1\wedge e_4\wedge e_5)\lrcorner d\psi&=&2{\rm Im}(\beta_{31})\\
 (e_2\wedge e_3\wedge e_5)\lrcorner d\psi&=&-2{\rm Re}(\beta_{12})\\
 (e_2\wedge e_4\wedge e_5)\lrcorner d\psi&=&2{\rm Im}(\beta_{32})\\
 (e_3\wedge e_4\wedge e_5)\lrcorner d\psi&=&-\mu+{\rm Im}(\beta_{11}+\beta_{22}-\beta_{33})\\
 \end{array}
 \right\}.
 \end{equation}
 Note that \[(e_2\wedge e_4)\lrcorner d\Omega-(e_1\wedge e_5)\lrcorner d\Omega-(e_2\wedge e_3\wedge e_5)\lrcorner d\psi-(e_1\wedge e_3\wedge e_4)\lrcorner d\psi=0.\] No other relations exist among the forms in (\ref{polarE3}), (\ref{polarE4}) and (\ref{polarE5}). Thus $s_5=9$ and $s_6=9$. Since $6s_0+5s_1+4s_2+3s_3+2s_4+s_5=34$, Cartan's test is satisfied and the proof is complete. 
 \end{proof}
 \begin{remark}
 The local nearly Calabi-Yau structures depend on $3$ functions of $6$ variables modulo diffeomorphisms.
 \end{remark}

 \section{Examples from twistor spaces of Riemannian four-manifolds}
 There has been an extensive literature on twistor theory. Suppose $(M^4, ds^2)$ is a Riemannian 4-manifold. A twistor at $x\in M$ is an orthogonal complex structure $\mathtt{j}:T_xM\rightarrow T_xM$, $\mathtt{j}^2=-1$ and $\mathtt{j}^*(ds_x^2)=ds_x^2.$ The space of twistors at points of $M$ forms a smooth manifold $\mathcal{J}$ called {\it twistor space} of $M$. It is well-known that $\mathcal{J}$ has an almost complex structure. It is moreover complex if $M$ has constant sectional curvatures. However, we will not use this usual almost complex structure in this paper. Instead, we will ``reverse" the almost complex structure on the fibers and obtain an $SU(3)$-structure on $\mathcal{J}$. By doing so, we will lose the possible integrability of the almost complex structures in some cases. 
 \subsection{Four dimensional Riemannian geometry}
 We formulate Riemannian geometry of four-manifolds in moving frames. Let $\pi:\mathcal{F}\rightarrow M$ be the oriented orthonormal coframe bundle over $M$. Thus $\mathcal{F}_x$ consists of orientation preserving isometries $u:T_xM\rightarrow \mathbf{R}^4$.
 Let $\eta$ be the $\mathbf{R}^4$-valued taugological form on $\mathcal{F}$. By the fundamental theorem of Riemannian geometry, there exists a unique $\mathfrak{so}(4,\mathbf{R})$-valued one-form $\theta$ so that 
 \[d\eta=-\theta\wedge\eta.\] 
 Denote $\omega_1=\eta_1+\sqrt{-1}\eta_3$ and $\omega_2=\eta_2+\sqrt{-1}\eta_4.$ We write the structrure equation as
 \begin{equation}\label{strequation4M}
 d\left(\begin{array}{c}
    \omega_1\\
    \omega_2\\
    \overline{\omega_1}\\
    \overline{\omega_2}
    \end{array}\right )=-
    \left(\begin{array}{cc}
          \alpha_{i\overline{j}}&\overline{\beta_{\overline{ij}}}\\
          \beta_{\overline{ij}}&\overline{\alpha_{i\overline{j}}}
          \end{array}\right)\wedge
     \left(\begin{array}{c}
    \omega_1\\
    \omega_2\\
    \overline{\omega_1}\\
    \overline{\omega_2}
    \end{array}\right )
 \end{equation}
 where $\alpha^t+\overline{\alpha}=0$ and $\beta^t+\beta=0$. 
 The Riemannian curvature is of course 
 \[R=d\left(\begin{array}{cc}
          \alpha &\overline{\beta}\\
          \beta &\overline{\alpha}
          \end{array}\right)+\left(\begin{array}{cc}
          \alpha &\overline{\beta}\\
          \beta &\overline{\alpha}
          \end{array}\right)\wedge\left(\begin{array}{cc}
          \alpha &\overline{\beta}\\
          \beta &\overline{\alpha}
          \end{array}\right).\]
   This is a $\mathfrak{so}(4,\mathbf{R})$-valued two form. Corresponding to the decomposition 
   $\mathfrak{so}(4,\mathbf{R})=\mathfrak{su}(2)_+\oplus\mathfrak{su}(2)_-$ we decompose $R=R_++R_-$, where
   \[R_+=d\left(\begin{array}{cc}
          \alpha_0 &0\\
           0&\overline{\alpha_0}
          \end{array}\right)+\left(\begin{array}{cc}
          \alpha_0 &0\\
          0&\overline{\alpha_0}
          \end{array}\right)\wedge\left(\begin{array}{cc}
          \alpha_0 &0\\
          0&\overline{\alpha_0}
          \end{array}\right)\] 
     and 
       \[R_-= d\left(\begin{array}{cc}
          \frac{1}{2}tr(\alpha)I &\overline{\beta}\\
          \beta &\frac{1}{2}\overline{tr(\alpha)I}
          \end{array}\right)+\left(\begin{array}{cc}
          \frac{1}{2}tr(\alpha)I &\overline{\beta}\\
          \beta &\frac{1}{2}\overline{tr(\alpha)I}
          \end{array}\right)\wedge\left(\begin{array}{cc}
          \frac{1}{2}tr(\alpha)I &\overline{\beta}\\
          \beta &\frac{1}{2}\overline{tr(\alpha)}I
          \end{array}\right).\]   
          for which $\alpha_0=\alpha-\frac{1}{2}tr(\alpha)I$ takes value in $\mathfrak{su}(2)$.  We are mainly interested in $R_-$ so we examine this part more carefully. Write 
          \[\beta=\left(\begin{array}{cc}
                       0&\omega_3\\
                       -\omega_3&0
                       \end{array}\right),\]
         \[(R_-)_1=\frac{1}{2}dtr(\alpha)+\omega_3\wedge\overline{\omega_3},\]
         and 
         \[(R_-)_2=d\omega_3-tr(\alpha)\wedge\omega_3.\]
   By the linear algebra developed in Section \ref{linearalgebra}, the forms \[\Theta_1=\omega_1\wedge\omega_2,\quad \Theta_2=\overline{\omega_1\wedge\omega_2},\quad  \Theta_3=\frac{\sqrt{-1}}{2}(\omega_1\wedge\overline{\omega_1}+\omega_2\wedge\overline{\omega_2})\]  form a basis for anti-self dual complex forms at $x$, while \[\Sigma_1=\omega_1\wedge\overline{\omega_2},\quad \Sigma_2=\overline{\omega_1}\wedge\omega_2,\quad  \Sigma_3=\frac{\sqrt{-1}}{2}(\omega_1\wedge\overline{\omega_1}-\omega_2\wedge\overline{\omega_2})\] form a basis for self dual forms at $x$. Since $R_-$ is semibasic, 
   \[(R_-)_1=A\Theta_1-\overline{A}\Theta_2+\sqrt{-1}a\Theta_3+B\Sigma_1-\overline{B}\Sigma_2+\sqrt{-1}b\Sigma_3,\]
   and 
   \[(R_-)_2=C_1\Theta_1+C_2\Theta_2+C_3\Theta_3+D_1\Sigma_1+D_2\Sigma_2+D_3\Sigma_3\]
   where $A, B, C_i, D_i$ are complex and $a, b$ are real.
  For our purposes, we view $R_-$ as a $(2,2)$ tensor. Using $ds^2$ and the convention of Section \ref{linearalgebra}, we write $R$ as 
  \[\begin{array}{ccl} 
  R_-&=&2(R_-)_1\otimes (E_{\bar{1}}\wedge\overline{E_{\bar{1}}}+ E_{\bar{2}}\wedge\overline{E_{\bar{2}}})+2(R_-)_2\otimes \overline{E_{\bar{1}}\wedge E_{\bar{2}}}+2\overline{(R_-)}_2\otimes E_{\bar{1}}\wedge E_{\bar{2}}\\\\
  &=&2\sqrt{-1}(R_-)_1\otimes\Theta_3^*+2(R_-)_2\otimes \Theta_2^*+2\overline{(R_-)}_2\otimes \Theta_1^*
  \end{array}\]
where, by abuse of notation, we use $E_{\bar{i}}$ to denote the tangent vector dual to $\omega_i$.
In this way we may regard $R_-$ as a linear map $R_-:\Lambda^2_-\rightarrow \Lambda^2=\Lambda^2_-\oplus \Lambda^2_+$. Relative to the basis $\Theta$ and $\Sigma$ we write the matrix representation 
\[R_-(\Theta_1,\Theta_2,\Theta_3)=2(\Theta_1,\Theta_2,\Theta_3,\Sigma_1,\Sigma_2,\Sigma_3)\left(
\begin{array}{ccc}
\overline{C_2}& C_1 &\sqrt{-1}A\\
\overline{C_1}&C_2  &-\sqrt{-1}\overline{A}\\ 
\overline{C_3}& C_3 &-a\\
\overline{D_2}& D_1&\sqrt{-1}B\\
\overline{D_1}&D_2&-\sqrt{-1}\overline{B}\\
\overline{D_3}&D_3&-b\\
\end{array}\right).
\]
 It is well-known that $R_-$ decomposes as $Z+W^-+\frac{s}{12}Id$ (see \cite{Besse}, p. 51) where $Z$ is the traceless Ricci curvature, $W^-$ is the anti-self-dual part of the Weyl curvature and $s$ is the scalar curvature. In our notations
 $Z$ is represented by the matrix 
 \[2\left(\begin{array}{ccc}
 \overline{D_2}& D_1&\sqrt{-1}B\\
\overline{D_1}&D_2&-\sqrt{-1}\overline{B}\\
\overline{D_3}&D_3&-b\\
\end{array}\right),
 \]  
  \[s=8(C_2+\overline{C_2}-a),\]
  and $W^-$  is represented by
  \[   2\left(
\begin{array}{ccc}
\overline{C_2}& C_1 &\sqrt{-1}A\\
\overline{C_1}&C_2  &-\sqrt{-1}\overline{A}\\ 
\overline{C_3}& C_3 &-a
\end{array}\right)-\frac{2}{3}(C_2+\overline{C_2}-a)I.\]

The metric with $W^-=0$ is called {\it self-dual}. If in addition, $ds^2$ is Einstein, then $s$ is necessarily constant. In our notations,
\begin{proposition}
The metric $ds^2$ is self-dual and Einstein if and only if $b=A=B=C_1=C_3=D_1=D_2=D_3=0$ and $C_2=\overline{C_2}=-a=\frac{s}{24}$. In this case, a part of the structure equation simplifies greatly 
 \begin{equation}\label{twistorstrequation}
 d\left(\begin{array}{c}
       \omega_1\\
       \omega_2\\
       \omega_3
       \end{array}\right)
 =-\left(\begin{array}{cc}
  \alpha&0 \\
  0&-tr(\alpha)\\
  \end{array}\right)
  \wedge\left(\begin{array}{c}
       \omega_1\\
       \omega_2\\
       \omega_3
       \end{array}\right)+
      \left(\begin{array}{c}
       \overline{\omega_2\wedge\omega_3}\\
       \overline{\omega_3\wedge\omega_1}\\
       \frac{s}{24}\overline{\omega_1\wedge\omega_2}
       \end{array}\right). \end{equation}
\end{proposition} 
Self-dual Einstein metrics will play important roles in our following constructions. There are not many compact examples with  $s\geq 0$ due to the classification by Hitchin (see \cite{Besse}, p 376):
\begin{theorem}
Let $M$ be be self-dual Einstein manifold. Then

(1) If $s>0$, $M$ is isometric to $S^4$ or $\mathbf{C}P^2$ with their canonical metrics.

(2) If $s=0$, $M$ is either flat or its universal covering is a $K3$ surface with the Calabi-Yau metric.
\end{theorem} 
The proof uses Bochner Technique, which, however does not work well when $s<0$. No similar results are available for self-dual Einstein metrics with negative scalar curvatures.
\subsection{Twistor spaces of self-dual Einstein manifolds}
We fix a complex structure $J_0$ on $\mathbf{R}^4$ by requiring $dz_1=dx_1+\sqrt{-1}dx_3$ and $dz_2=dx_2+\sqrt{-1}dx_4$ be complex linear. We define a map $\mathtt{j}:\mathcal{F}\rightarrow \mathcal{J}$ as follows 
\[\mathtt{j}(u)=u^{-1}\circ J_0\circ u.\] 
Since $SO(4)$ acts transitively on the orthogonal complex structures on $\mathbf{R}^4$ and the isotropic group of $J_0$ is $U(2)$, $\mathtt{j}$ makes $\mathcal{F}$ a principal $U(2)$-bundle over $\mathcal{J}$. This defines a $U(2)$-structure on $\mathcal{J}$. It in turn determines an $SU(3)$-structure on $\mathcal{J}$ by the standard embedding of $U(2)$ into $SU(3)$. Relative to $\mathtt{j}$, $\beta$ becomes semi-basic. The almost structure on $\mathcal{J}$ determined by this $SU(3)$-strucure is such that $\omega_1, \omega_2$ and $\omega_3$ are complex linear. 

Let us now concentrate on self-dual Einstein manifolds. The equations in (\ref{twistorstrequation}) are the first structure equations on $\mathcal{J}$. It clearly satisfies the condition of Lemma \ref{involutivity}. Thus special Lagrangian geometry behaves nicely on the twistor spaces of these manifolds.
\subsubsection{$s>0$} 
In this case, we scale the metric so that $s=24$. Now the structure equation (\ref{twistorstrequation}) indicates that the twistor space is actually nearly K\"{a}hler. By the aforementioned Hitchin's result, the only two possibilities are $M=S^4$ and $M=\mathbf{CP}^2$. The corresponding twistor spaces are two familiar nearly K\"{a}hler examples, $\mathbf{CP}^3$ and the flag manifold $SU(3)/T^2$.  
\subsubsection{$s=0$} Again, by Hitchin's result we have two examples: one is the flat case, the other is $K_3$ surfaces. 
\subsubsection{$s<0$} This is the most interesting case in many aspects. We scale the metric to make $s=-48$. Now the structure equation (\ref{twistorstrequation}) reads
\begin{equation}
 d\left(\begin{array}{c}
       \omega_1\\
       \omega_2\\
       \omega_3
       \end{array}\right)
 =-\left(\begin{array}{cc}
  \alpha&0 \\
  0&-tr(\alpha)\\
  \end{array}\right)
  \wedge\left(\begin{array}{c}
       \omega_1\\
       \omega_2\\
       \omega_3
       \end{array}\right)+
      \left(\begin{array}{c}
       \overline{\omega_2\wedge\omega_3}\\
       \overline{\omega_3\wedge\omega_1}\\
       -2\overline{\omega_1\wedge\omega_2}
       \end{array}\right). \end{equation}
   The only torsion is the Nijenhuis tensor, in local unitary basis,  
   \[N=\left(\begin{array}{ccc}
       1&0&0\\
        0&1&0\\
        0&0&-2
   \end{array}\right). 
     \]    
Thus, we have 
\begin{theorem}
The twistor space of a self-dual Einstein manifold of negative scalar curvature is strictly nearly Calabi-Yau. 
\end{theorem}
The simplest example of this category is, of course, the twistor space of the hyperbolic space $H^4$. Compact examples can be obtained from the quotients of $H^4$ by certain discrete isometry groups.  
\section{Complete special Lagrangian examples}
In this section we will construct some complete special Lagrangian submanifolds in the twistor spaces $\mathcal{J}(S^4)=\mathbf{CP}^3$ and $\mathcal{J}(H^4)$ considered in the previous section. Our method is based on the following observation due to Robert Bryant \cite{Bryant2}. Suppose on an $SU(3)$ manifold $(M,\Omega,\Psi)$, there is a real structure, i.e., an involution $c$ such that 
\[c^*\Omega=-\Omega,\quad  c^*\Psi=\overline{\Psi},\]
and the set $N_c$ of points fixed under $c$ is a smooth submanifold. Then it is easy to see that $N_c$, with one of its two possible orientations, is a special Lagrangian submanifold of $M$ . Thus our major task is to construct such involutions for $\mathcal{J}(S^4)=\mathbf{CP}^3$ and $\mathcal{J}(H^4)$.
\subsection{An example in $\mathcal{J}(S^4)=\mathbf{CP}^3$}
We need a more explicit description of the twistor fibration $T:\mathbf{CP}^3\rightarrow S^4$. We follow the discussion in \cite{Bryant3}. However, as aformentioned, we will use a different almost complex structure on $\mathbf{CP}^3$.

Let $\mathbf{H}$ denote the real division algebra of quaternions. An element of $\mathbf{H}$ can be written uniquely as $q=z+jw$ where $z,w\in \mathbf{C}$ and $j\in\mathbf{H}$ satisfies 
  \[j^2=-1,\quad  zj=j\bar{z}\]
  for all $z\in\mathbf{C}.$ The quaternion multiplication is thus given by
  \begin{equation}\label{productrule}(z_1+jz_2)(z_3+jz_4)=z_1z_3-z_2\overline{z_4}+j(z_2z_3+\overline{z_1z_4}).\end{equation}We define an involution $C:\mathbf{H}\rightarrow \mathbf{H}$ by $C(z_1+jz_2)=\overline{z_1}+j\overline{z_2}$. It can be easily checked via the product rule (\ref{productrule}) that this is in fact an algebra automorphism, i.e., $C(p q)=C(p)C(q)$ for $p,q \in \mathbf{H}$. 

We regard $\mathbf{C}$ as subalgebra of $\mathbf{H}$ and give $\mathbf{H}$ the structure of a complex vector space by letting $\mathbf{C}$ act on the right. 
  We let $\mathbf{H}^2$ denote the space of pairs $(q_1,q_2)$ where $q_i\in\mathbf{H}$. We will make $\mathbf{H}^2$ into a quaternion vector space by letting $\mathbf{H}$ act on the right 
  \[(q_1,q_2)q=(q_1q,q_2q).\]
This automatically makes $\mathbf{H}^2$ into a complex vector space of dimension $4$. In fact, regarding $\mathbf{C}^4$ as the space of 4-tuples $(z_1,z_2,z_3,z_4)$, we make the explicit identification 
 \begin{equation}\label{H=C^2}(z_1,z_2,z_3,z_4)\sim(z_1+jz_2,z_3+jz_4).\end{equation}
 This specific isomorphism is the one we will always mean when we write $\mathbf{C}^4=\mathbf{H}^2$. 

 If $v\in\mathbf{H}^2\setminus(0,0)$ is given, let $v\mathbf{C}$ and $v\mathbf{H}$ denote, respectively, the complex line and the quaternion line spanned by $v$. As is well-known, $\mathbf{HP}^1$, the space of quaternion lines in $\mathbf{H}^2$, is isometric to $S^4.$ For this reason, we will speak interchangeably of $S^4$ and $\mathbf{HP}^1$.  The assignment $v\mathbf{C}\rightarrow v\mathbf{H}$ is exactly the twistor mapping $T: \mathbf{CP}^3\rightarrow\mathbf{HP}^1$. The fibres of $T$ are $\mathbf{CP}^1$'s. Thus, we have a fibration 
 \begin{equation}\label{fibration1}
 \begin{array}{lcc}
    \mathbf{CP}^1&\rightarrow &\mathbf{CP}^3\\
                 &            &\downarrow\\
                 &            &\mathbf{HP}^1\\ 
                 \end{array}\end{equation}
This is the famous twistor fibration. In order to study its geometry more thoroughly, we will now introduce the structure equations of $\mathbf{H}^2$.
First we endow $\mathbf{H}^2$ with a quaternion inner product $\left\langle, \right\rangle: \mathbf{H}^2\times\mathbf{H}^2\rightarrow \mathbf{H}$ defined by 
 \[\left\langle (q_1,q_2), (p_1,p_2)\right\rangle=\bar{q_1}p_1+\bar{q_2}p_2.\]
 We have identities
 \[\left\langle v,wq\right\rangle =\left\langle v,w\right\rangle q,\quad \overline{\left\langle v,w\right\rangle}=\left\langle w,v\right\rangle, \quad \left\langle vq,w\right\rangle=\bar{q}\left\langle v,w\right\rangle.\]
Via the identification (\ref{H=C^2}), \begin{equation}\label{qproduct}\begin{array}{ccl}
\langle(q_1,q_2),(p_1,p_2)\rangle&=&\overline{z_1}w_1+\overline{z_2}w_2+\overline{z_3}w_3+\overline{z_4}w_4\\
&&+j(z_1w_2-z_2w_1+z_3w_4-z_4w_3)              \end{array}          \end{equation}
for $q_1=z_1+jz_2,\quad q_2=z_3+jz_4,\quad p_1=w_1+jw_2,\quad p_2=w_3+jw_4.$ In other words, $\langle,\rangle$ essentially consists of two parts: one is the standard Hermitian product $dz_1\circ d\overline{z_1}+\cdots+dz_4\circ d\overline{z_4}$; the other is the standard complex symplectic form $dz_1\wedge dz_2+dz_3\wedge dz_4$. 
 
 Let $\mathfrak{F}$ denote the space of pairs $f=(e_1,e_2)$ with $e_i\in \mathbf{H}^2$ satisfying 
 \[\left\langle e_1,e_1 \right\rangle=\left\langle e_2, e_2\right\rangle=1, \quad \left\langle e_1,e_2\right\rangle=0.\]
 We regard $e_i(f)$ as functions on $\mathfrak{F}$ with values in $\mathbf{H}^2$. Clearly $e_1(\mathfrak{F})=S^7\subset \mathbf{E}^8=\mathbf{H}^2$.
 It is well-known that $\mathfrak{F}$ may be canonically identified with $Sp(2)$ up to a left translation in $Sp(2)$. There are unique quaternion-valued 1-forms $\{\phi^{a}_{b}\}$ so that 
 \begin{equation}\label{strequation1}
 de_a=e_b\phi^b_{a},
 \end{equation}
 
  \begin{equation}\label{strequation2}
  d\phi^a_{b}+\phi^a_c\wedge\phi^c_b=0,
  \end{equation} 
    and 
    \begin{equation}
  \phi^a_{b}+\overline{\phi^b_a}=0.
  \end{equation}

We have a canonical map $\mathfrak{F}\rightarrow \mathbf{CP}^3$ by sending $(e_1,e_2)$ to the complex line spanned by $e_1$. We will denote this map by $\mathtt{j}$ by a slight abuse of notation. The composition $\pi=T\circ \mathtt{j}$ is actually a spin structure on $S^4$. In fact the oriented coframe bundle $\mathbf{F}$ of $S^4$ may be identified with $SO(5)$ up to a left translation in $SO(5)$. Thus $\mathfrak{F}$ double covers $\mathbf{F}$ as $Sp(2)$ double covers $SO(5)$. 

We now write structure equations for the map $\mathtt{j}$. First we immediately see that $\mathtt{j}$ gives $\mathfrak{F}$ an $S^1\times S^3$-structure over $\mathbf{CP}^3$ where we have identified $S^1$ with the unit complex numbers and $S^3$ with the unit quaternions. The action is given by 
\[f(z,q)=(e_1,e_2)(z,q)=(e_1z,e_2q)\] where $z\in S^1$ and $q\in S^3$. If we set 
\[\left[\begin{array}{ll} 
     \phi^1_{1}&\phi^1_2\\
     \phi^2_1 &\phi^2_2\\
     \end{array}\right]=
   \left[\begin{array}{cc}
     i\rho_1+j\overline{\omega_3}& -\frac{\overline{\omega_1}}{\sqrt{2}}+j\frac{\omega_2}{\sqrt{2}}\\
     \frac{\omega_1}{\sqrt{2}}+j\frac{\omega_2}{\sqrt{2}}& i\rho_2+j\tau\\
     \end{array}\right]   
   \]
where $\rho_1$ and $\rho_2$ are real 1-forms while $\omega_1$, $\omega_2$, $\omega_3$ and $\tau$ are complex valued,  we may rewrite one part of the structure equation (\ref{strequation2}) relative to the $S^1\times S^3$-structure on $\mathbf{CP}^3$ as
   \begin{equation}\label{strequation4}d\left(
   \begin{array}{l}
   \omega_1\\
    \omega_2\\
    \omega_3\\
    \end{array}\right)=-\left(
    \begin{array}{lcl}
    i(\rho_2-\rho_1)&-\bar{\tau}&0\\
    \tau& -i(\rho_1+\rho_2)&0\\
    0&0&2i\rho_1\\
    \end{array} \right)\wedge\left(
    \begin{array}{l}
    \omega_1\\\omega_2\\\omega_3\end{array}\right)+\left(
    \begin{array}{l}    \overline{\omega_2\wedge\omega_3}\\\overline{\omega_3\wedge\omega_1}\\\overline{\omega_1\wedge\omega_2}\end{array}\right).
 \end{equation}  
The nearly K\"{a}hler structure on $\mathbf{CP}^3$ is defined by setting $\omega_1$, $\omega_2$ and $\omega_3$ to be complex linear. 

Via the algebra automorphism $C$ we define an involution on $\mathbf{H}^2$ by $(p,q)\mapsto (C(p),C(q))$. We denote this map still by $C$. This map in turn induces an involution on $\mathfrak{F}$, still denoted $C$, by $C(e_1,e_2)=(C(e_1),C(e_2))$. From (\ref{qproduct}) we see that the defining equations for $\mathfrak{F}$ are preserved and the involution is well-defined. The map $C$ further descends to an involution $c$ on $\mathbf{CP}^3$ as well as an involution $\bar{c}$ on $S^4$ by  $e_1\mathbf{C}\mapsto C(e_1)\mathbf{C}$ and $e_1\mathbf{H}\mapsto C(e_1)\mathbf{H}$ repectively. We have the following commutative diagram
\[\begin{array}{ccc}
   \mathfrak{F} &\stackrel{C}{\longrightarrow}&\mathfrak{F}\\
    \downarrow&&\downarrow\\
    \mathbf{CP}^3&\stackrel{c}{\longrightarrow}&\mathbf{CP}^3\\
\downarrow&&\downarrow\\
    \mathbf{HP}^1&\stackrel{\bar{c}}{\longrightarrow}&\mathbf{HP}^1\\
  \end{array}.
\]
Apply the automorphism to the structure equations (\ref{strequation1}) and we get
\[dC(e_a)=C(e_b)C(\phi_{a}^b).\]
Thus in particular we have on $\mathfrak{F}$
\[C^*\omega_i=\overline{\omega_i}\] for $i=1,2,3$.
Consequently 
 \[C^*\mathtt{j}^*\Omega=-\mathtt{j}^*\Omega, \quad C^*\mathtt{j}^*\Psi=\overline{\mathtt{j}^*\Psi}. \]
Since $\mathtt{j}C=c\mathtt{j}$ and $\mathtt{j}^*$ is injective, we have on $\mathbf{CP}^3$
\[c^*\Omega=-\Omega, \quad c^*\Psi=\overline{\Psi}.\]

Thus by the general principle the fixed set of $c$ is a special Lagrangian submanifold of $\mathbf{CP}^3$. Moreover it is easy to see that this locus is just the usual $\mathbf{RP}^3$.
\begin{theorem}
The real projective space $\mathbf{RP}^3=\{[x_1:x_2:x_3:x_4]:x_i\in\mathbf{R}\}\subset\mathbf{CP}^3$ is a special Lagrangian submanifold of the nearly K\"{a}hler $\mathbf{CP}^3$. 
\end{theorem}
 The twistor map $T:\mathbf{CP}^3\rightarrow\mathbf{HP}^1$ restricted to the real projective space $\mathbf{RP}^3$ now looks like
  \[[x_1:x_2:x_3:x_4]\mapsto [x_1+jx_2:x_3+jx_4].\]
Thus the image is a $\mathbf{CP}^1\subset\mathbf{HP}^1$ and $T$ is the Hopf fibration 
    \[\begin{array}{ccc}
       S^1&\rightarrow&\mathbf{RP}^3\\
          &&            \downarrow\\
          &&\mathbf{CP}^1
      \end{array}.
\]

A dual construction for $\mathcal{J}(H^4)$ will follow.
\subsection{An example in $\mathcal{J}(H^4)$}
Let $\mathbf{H}^2$ and the involution $C$ be as before. But now we endow $\mathbf{H}^2$ with a $(1,1)$ quaternion inner product $\langle,\rangle:\mathbf{H}^2\times \mathbf{H}^2\rightarrow \mathbf{H}$ defined by 
\[\langle(q_1,q_2),(p_1,p_2)\rangle=\overline{q_1}p_1-\overline{q_2}p_2.\]
We still have the identities
\[\left\langle v,wq\right\rangle =\left\langle v,w\right\rangle q, \quad \overline{\left\langle v,w\right\rangle}=\left\langle w,v\right\rangle, \quad \left\langle vq,w\right\rangle=\bar{q}\left\langle v,w\right\rangle.\]
Via the identification (\ref{H=C^2}), \begin{equation}\label{pqproduct}\begin{array}{ccl}
\langle(q_1,q_2),(p_1,p_2)\rangle&=&\overline{z_1}w_1+\overline{z_2}w_2-\overline{z_3}w_3-\overline{z_4}w_4\\
&&+j(z_1w_2+z_2w_1-z_3w_4-z_4w_3)              \end{array}          \end{equation}
for $q_1=z_1+jz_2,\quad q_2=z_3+jz_4,\quad p_1=w_1+jw_2,\quad p_2=w_3+jw_4.$ In other words, $\langle,\rangle$ essentially consists of two parts: one is the $(2,2)$ Hermitian product $dz_1\circ d\overline{z_1}+dz_2\circ d\overline{z_2}-dz_3\circ d\overline{z_3}-dz_4\circ d\overline{z_4}$; the other is a complex symplectic form $dz_1\wedge dz_2-dz_3\wedge dz_4$.
Denote the pseudo-sphere in $\mathbf{H}^2$ by
  \[\Psi S^7=\{(p,q)\in \mathbf{H}^2: \overline{p}p-\overline{q}q=1\}.\] This is a connected non-compact smooth hypersurface in $\mathbf{H}^2$. The group $S^3$ acts on $\Psi S^7$ by \[(p,q)\cdot r=(pr,qr)\] 
where $r\in S^3$ is a unit quaternion number. This action is clearly free. Thus the quotient space $\Psi \mathbf{HP}^1=\Psi S^7/S^3$ is smooth. Indeed $\Psi\mathbf{HP}^1=H^4$. Similarly if we regard $S^1$ as a subgroup of $S^3$ consisting of unit complex numbers, the quotient space $\Psi \mathbf{CP}^3=\Psi S^7/S^1$ is smooth. The clearly well-defined map $T:\Psi\mathbf{CP}^3\rightarrow H^4$ is exactly the twistor fibration of $H^4$. We have the following commutative diagram of fibrations
                \[\begin{array}{ccc}
                  S^1&\hookrightarrow& \Psi S^7\\
                  &&\downarrow\\ 
                  \mathbf{CP}^1&\hookrightarrow&\Psi\mathbf{CP}^3\\
    &&\downarrow\\                
&&\Psi\mathbf{HP}^1      
                 \end{array}.
\] 

 Let $\mathfrak{F}$ denote the space of pairs $f=(e_1,e_2)$ with $e_i\in \mathbf{H}^2$ satisfying 
 \[\left\langle e_1,e_1 \right\rangle=1, \quad \left\langle e_2, e_2\right\rangle=-1, \quad \left\langle e_1,e_2\right\rangle=0.\]
 We regard $e_i(f)$ as functions on $\mathfrak{F}$ with values in $\mathbf{H}^2$. Clearly $e_1(\mathfrak{F})=\Psi S^7\subset $E$^{(4,4)}=\mathbf{H}^2$.
 It is well-known that $\mathfrak{F}$ maybe canonically identified with $Sp(1,1)$ up to a left translation in $Sp(1,1)$, where 
\[Sp(1,1)=\{A\in \mathfrak{gl}(2,\mathbf{H}): \overline{A}\left(\begin{array}{cc}
                                                                  1 &0\\
                                                                  0&-1\\
                                                                \end{array}
\right)A^t=\left(\begin{array}{cc}
                                                                  1 &0\\
                                                                  0&-1\\
                                                                \end{array}\right).\}\]
 There are unique quaternion-valued 1-forms $\{\phi^{a}_{b}\}$ so that 
 \begin{equation}\label{pstrequation1}
 de_a=e_b\phi^b_{a},
 \end{equation}
 
  \begin{equation}\label{pstrequation2}
  d\phi^a_{b}+\phi^a_c\wedge\phi^c_b=0,
  \end{equation} 
    and 
    \begin{equation}
  \overline{\phi}\left(\begin{array}{cc}
                                                                  1 &0\\
                                                                  0&-1\\
                                                                \end{array}
\right)+\left(\begin{array}{cc}
                                                                  1 &0\\
                                                                  0&-1\\
                                                                \end{array}
\right)\phi^t=0.
  \end{equation}

We have a canonical map $\mathfrak{F}\rightarrow \Psi\mathbf{CP}^3$ by sending $(e_1,e_2)$ to the coset $e_1\cdot S^1$. We will denote this map by $\mathtt{j}$ by a slight abuse of notation. The composition $\pi=T\circ \mathtt{j}$ is actually a spin structure on $H^4$. In fact the oriented coframe bundle $\mathbf{F}$ of $H^4$ may be identified with $SO^0(4,1)$, the identity component of $SO(4,1)$, up to a left translation in $SO^0(4,1)$. Thus $\mathfrak{F}$ double covers $\mathbf{F}$ as $Sp(1,1)$ double covers $SO^0(4,1)$ (see Harvey \cite{Harvey}, p. 272 for the isomorphism $Sp(1,1)\cong Spin^0(4,1)$ where he used the notation $HU(1,1)$ for $Sp(1,1)$). 

We now write structure equations for the map $\mathtt{j}$. First we immediately see that $\mathtt{j}$ gives $\mathfrak{F}$ an $S^1\times S^3$-structure over $\Psi\mathbf{CP}^3$ where we have identified $S^1$ with the unit complex numbers and $S^3$ with the unit quaternions. The action is given by 
\[f(z,q)=(e_1,e_2)\cdot (z,q)=(e_1z,e_2q)\] where $z\in S^1$ and $q\in S^3$. If we set 
\[\left[\begin{array}{ll} 
     \phi^1_{1}&\phi^1_2\\
     \phi^2_1 &\phi^2_2\\
     \end{array}\right]=
   \left[\begin{array}{cc}
     i\rho_1+j\overline{\omega_3}& \overline{\omega_1}-j\omega_2\\
     \omega_1+j\omega_2& i\rho_2+j\tau\\
     \end{array}\right]   
   \]
where $\rho_1$ and $\rho_2$ are real 1-forms while $\omega_1$, $\omega_2$, $\omega_3$ and $\tau$ are complex valued,  we may rewrite one part of the structure equation (\ref{strequation2}) relative to the nearly Calabi-Yau structure on $\Psi\mathbf{CP}^3$ as
   \begin{equation}\label{pstrequation4}d\left(
   \begin{array}{l}
   \omega_1\\
    \omega_2\\
    \omega_3\\
    \end{array}\right)=-\left(
    \begin{array}{lcl}
    i(\rho_2-\rho_1)&-\bar{\tau}&0\\
    \tau& -i(\rho_1+\rho_2)&0\\
    0&0&2i\rho_1\\
    \end{array} \right)\wedge\left(
    \begin{array}{r}
    \omega_1\\\omega_2\\\omega_3\end{array}\right)+\left(
    \begin{array}{r}    \overline{\omega_2\wedge\omega_3}\\\overline{\omega_3\wedge\omega_1}\\-2\overline{\omega_1\wedge\omega_2}\end{array}\right).
 \end{equation}  
The nearly Calabi-Yau structure on $\Psi\mathbf{CP}^3$ is defined by setting $\omega_1$, $\omega_2$ and $\omega_3$ to be complex linear. 

 This involution $C$ on $\mathbf{H}^2$ induces an involution on $\mathfrak{F}$, still denoted by $C$, by $C(e_1,e_2)=(C(e_1),C(e_2))$. From (\ref{pqproduct}) we see that the defining equations for $\mathfrak{F}$ are preserved and the involution is well-defined. The map $C$ further descends to an involution $c$ on $\Psi\mathbf{CP}^3$ as well as an involution $\bar{c}$ on $H^4$ by $e_1\cdot S^1\mapsto C(e_1)\cdot S^1$ and $e_1\cdot S^3\mapsto C(e_1)\cdot S^3$ repectively. We have the following commutative diagram
\[\begin{array}{ccc}
   \mathfrak{F} &\stackrel{C}{\longrightarrow}&\mathfrak{F}\\
    \downarrow&&\downarrow\\
    \Psi\mathbf{CP}^3&\stackrel{c}{\longrightarrow}&\Psi\mathbf{CP}^3\\
\downarrow&&\downarrow\\
    H^4&\stackrel{\bar{c}}{\longrightarrow}&H^4\\
  \end{array}.
\]
Apply the automorphism to the structure equations (\ref{pstrequation1}) and we get
\[dC(e_a)=C(e_b)C(\phi_{a}^b).\]
Thus in particular we have on $\mathfrak{F}$
\[C^*\omega_i=\overline{\omega_i}\] for $i=1,2,3$.
Consequently 
 \[C^*\mathtt{j}^*\Omega=-\mathtt{j}^*\Omega,\quad C^*\mathtt{j}^*\Psi=\overline{\mathtt{j}^*\Psi}. \]
Since $\mathtt{j}C=c\mathtt{j}$ and $\mathtt{j}^*$ is injective, we have on $\Psi\mathbf{CP}^3$
\[c^*\Omega=-\Omega,\quad c^*\Psi=\overline{\Psi}.\]

Thus by the general principle the fixed set of $c$ is a special Lagrangian submanifold of $\Psi\mathbf{CP}^3$. It is easy to see that this manifold is the pseudo-projective $3$-space $\Psi\mathbf{RP}^3$, defined as the quotient of the pseudo 3-sphere $\Psi S^3=\{(x_1,x_2,x_3,x_4)\in \mathbf{R}^4: x_1^2+x_2^2-x_3^2-x_4^2=1\}$ by $\mathbb{Z}_2$. 

\begin{theorem}
The real pseudo-projective space $\Psi\mathbf{RP}^3\subset\Psi\mathbf{CP}^3$ is a special Lagrangian submanifold of the nearly Calabi-Yau $\Psi\mathbf{CP}^3$. 
\end{theorem}
The image under $T$ of this pseudo-sphere is easily seen to be the hyperbolic $2$-space $H^2\subset H^4.$ Thus we have the following fibration 
\[\begin{array}{ccc}
       S^1&\rightarrow&\Psi\mathbf{RP}^3\\
          &&            \downarrow\\
          &&H^2
      \end{array}.
\]

 \section{Compact special Lagrangian submanifolds in nearly Calabi-Yau manifolds} 
 We discuss compact special Lagrangian submainifolds in nearly Calabi-Yau manifolds. We answer two questions:
\begin{enumerate}
 \item Let $N$ be a compact special Lagrangian $3$-fold in a fixed nearly Calabi-Yau manifold $(M,\Omega,\Psi)$. Let $\mathcal{M}_N$ be the moduli space of special Lagrangian deformations of $N$, that is, the connected component of the set of special Lagrangian $3$-folds containing $N$. What can we say about $\mathcal{M}_N$? Is it a smooth manifold? what is its dimension?
 \item Let $\{(M, \Omega_t, \Psi_t: t\in (-\epsilon,\epsilon)\}$ be a smooth $1$-parameter family of nearly Calabi-Yau manifolds. Suppose $N_0$ is an SL-$3-$fold. Under what conditions can we extend $N_0$ to a smooth family of special Lagrangian $3$-folds $N_t$ in $(M, \Omega_t,\Psi_t)$?
 
 \end{enumerate}
 These questions concern the deformations of special Lagrangian $3$-folds and obstructions to their existence repecitively. In the Calabi-Yau case, the first question is answered by R. McLean in \cite{McLean}, and the second is answered by D. Joyce in \cite{Joyce}. Moreover, \cite{Joyce} also answers these questions for more general K\"{a}hler manifolds, which Joyce called {\it almost Calabi-Yau}. We show that their proofs generalize to nearly Calabi-Yau manifolds.
\subsection{Deformations of compact special Lagrangian $3$-folds}
We have the following result similar to \cite{McLean}
\begin{theorem}
Let $(M,\Omega,\Psi)$ be a nearly Calabi-Yau $3$-fold, and $N$ a compact special Lagrangian $3$-fold in $M$. Then the moduli space $\mathcal{M}_N$ of special Lagrangian deformations of $N$ is a smooth manifold of dimension $b^1(N)$, the first Betti number of $N$.
\end{theorem} 
The proof is a word-for-word copy of \cite{McLean}. Note that McLean's proof only depends on the fact that $\Omega$ and ${\rm Im}\Psi$ are both closed. These conditions are satisfied by nearly Calabi-Yau manifolds.
\subsection{Obstructions to the existence of compact SL $3$-folds}
We address Question 2 above. Let $\{(M,\Omega_t,\Psi_t)\}$ be a smooth $1$-parameter family of nearly Calabi-Yau manifolds. Supoose $N_0$ is a special Lagangian $3$-fold of $(M,\Omega_0,\Psi_0)$ and $N_t$ is an extension. Then we can view $N_t$ as a family of embeddings of $\mathbf{i}_t: N_0\rightarrow M$ such that $\mathbf{i}_t^*(\Omega_t)=\mathbf{i}_t^*({\rm Im}\Psi_t)=0$. Since the cohomology classes $[\mathbf{i}_s^*(\Omega_t)]$ and $[\mathbf{i}_s^*({\rm Im}\Psi_t)]$ do not vary with $s$, we have $[\mathbf{i}_0^*(\Omega_t)]=[\mathbf{i}_0^*({\rm Im}\Psi_t)]=0$. Thus a necessary condition for such an extension of $N_0$ to exist is 
\[[\Omega_t|_{N_0}]=[{\rm Im}\Psi_t|_{N_0}]=0.\] Actually this is also sufficient.
\begin{theorem}
Let  $\{(M,\Omega_t,\Psi_t): t\in (-\epsilon,\epsilon)\}$ be a smooth $1$-parameter family of nearly Calabi-Yau $3$-folds. Let $N_0$ be a compact SL $3$-fold in $(M,\Omega_0,\Psi_0)$, and suppose $[\Omega_t|_{N_0}]=0$ in $H^2(N_0,\mathbf{R})$ and $[{\rm Im}\Psi_t|_{N_0}]=0$ in $H^3(N_0,\mathbf{R})$ for all $t\in (-\epsilon, \epsilon)$. Then $N_0$ extends to a smooth $1$-parameter family $\{N_t:t\in (-\delta, \delta)\}$ for some $0<\delta\leq \epsilon$ and $N_t$ is a compact SL $3$-fold in $(M,\Omega_t,\Psi_t)$.
\end{theorem}
Again, the proof follows exactly the Calabi-Yau case as in \cite{Joyce}. However, since the details are not readily available, we write them down. 
\begin{proof} Let $\nu_{N_0}$ be the normal bundle of $N_0$ in $(M,\Omega_0,\Psi_0)$. Denote $\mathtt{exp}$ the exponential map of $(M,\Omega_0,\Psi_0)$. For a vector bundle $E$ over $N_0$, we use $C^{1,\alpha}(E)$ and $C^{0,\alpha}(E)$ to denote the sections of $E$ of class $C^{1,\alpha}$ and $C^{0,\alpha}$ respectively. We define a map 
\[ F:C^{1,\alpha}(\nu_{N_0})\times (-\epsilon,\epsilon)\rightarrow dC^{1,\alpha}(\Lambda^1(N_0))\times dC^{1,\alpha}(\Lambda^2(N_0))\] by
 \[F(V,t)=(\mathtt{exp}_V^*(\Omega_t),\mathtt{exp}_V^*({\rm Im}\Psi_t)).\]
We need show this is well-defined. The maps $F(sV,t)(0\leq s\leq 1)$ provide a homotopy between $F(0,t)=(\Omega_t|_{N_0}, {\rm Im}\Psi|_{N_0})$ and $F(V,t)$. Since $[\Omega_t|_{N_0}]=0$, $[\mathtt{exp}_V^*(\Omega_t)]=0$. Thus $\mathtt{exp}_V^*(\Omega_t)=d\tau$ for some $\tau$. Moreover, by the standard Hodge theory the form $\tau$ can be chosen to be in $C^{2,\alpha}$ because $V$ is $C^{1,\alpha}$ and so is $\mathtt{exp}_V^*(\Omega_t)$. It is similar to show that $\mathtt{exp}_V^*({\rm Im}\Psi_t)$ lies in $dC^{1,\alpha}(\Lambda^2(N_0))$. 

Now we compute the tangent map of $F$ at the point $(0,0)$, \[F_*:\mathbf{R}\times C^{1,\alpha}(\nu_{N_0})\rightarrow dC^{1,\alpha}(\Lambda^1(N_0))\times dC^{1,\alpha}(\Lambda^2(N_0)).\]
First\[\begin{array}{ccl}F_*(\frac{\partial}{\partial t},0)&=&\frac{d}{dt}|_{t=0,V=0}(\mathtt{exp}_V^*\Omega_t, \mathtt{exp}_V^*{\rm Im}\Psi_t)\\\\
        &=&(\dot{\Omega}|_{N_0},{\rm Im}(\dot{\Psi})|_{N_0}).
       \end{array}
\] where \[\dot{\Omega}=\frac{d}{dt}|_{t=0}\Omega_t,\quad \dot{\Psi}=\frac{d}{dt}|_{t=0}\Psi_t.\]
Second\[\begin{array}{ccl}
   F_*(0,V)&=&\frac{d}{ds}|_{s=0}(\mathtt{exp}_{sV}^*\Omega_0,\mathtt{exp}_{sV}^*\Omega_0)\\\\
           &=&(\mathcal{L}_V\Omega_0|_{N_0},\mathcal{L}_V{\rm Im}\Psi_0|_{N_0})\\\\
           &=&((V\lrcorner d\Omega_0+d(V\lrcorner\Omega_0))|_{N_0},(V\lrcorner d({\rm Im}\Psi_0)+d(V\lrcorner {\rm Im}\Psi_0))|_{N_0})\\\\
           &=&(d(V\lrcorner\Omega_0)|_{N_0},d(V\lrcorner {\rm Im}\Psi_0)|_{N_0})\\
\end{array}
\]where $\mathcal{L}_V$ is the Lie derivative in the $V$ direction and the Cartan formula is used. Actually, in order to take the Lie derivative, one must extend the normal vector field $V$ to an open neighborhood first. It is easy to see the result is independent of this extension. 

Note that the mapping $V\mapsto v=V\lrcorner \Omega_0$ gives a bundle isomporphism between $T^*N_0$ and $\nu_{N_0}$. Translated via this correspondence $V\lrcorner {\rm Im}\Psi$ becomes $-*v$ as is shown in \cite{McLean} where $*$ is the Hodge star operation. By Hodge theory, $F_*(0,V)$ runs over every element in $dC^{1,\alpha}(\Lambda^1(N_0))\times dC^{1,\alpha}(\Lambda^2)(N_0)$. Thus $F_*$ is surjective. 
We can also compute  the kernerl 
\[
F_*^{-1}(0,0)=\{(r\frac{\partial}{\partial t},V): r\dot{\Omega}|_{N_0}=-dv, r{\rm Im}(\dot{\Psi})|_{N_0}=d*v, r\in \mathbf{R}\},\] where $v$ relates to $V$ as above. Since $[\Omega_t|_{N_0}]=0$ and $[{\rm Im}(\Psi_t)|_{N_0}]=0$ we have $[\dot{\Omega}|_{N_0}]=0$ and $[{\rm Im}(\dot{\Psi})|_{N_0}]=0$. Thus $\dot{\Omega}|_{N_0}$ and ${\rm Im}(\dot{\Psi})|_{N_0}$ are exact. Again by Hodge theory, $F_{*}^{-1}(0,0)$ is nonempty and finite-dimensional with dimension $b^1(N_0)+1$. By the Implicit Function Theorem for smooth maps between Banach spaces, $F^{-1}(0,0)$ is a smooth manifold with its tangent space at $(0,0)$ equal to $F_{*}^{-1}(0,0)$. The $C^{1,\alpha}(\gamma_{N_0})$ components of elements of $F^{-1}(0,0)$ are in fact smooth sections by the elliptic regularity theorem. Note that the projection map $t$ restricted to $F^{-1}(0,0)$ is nondegenerate at $(0,0)$. Thus the manifold $F^{-1}(0,0)$ is a local smooth fibration over $(-\epsilon,\epsilon)$. Pick a local section $(t,V_t)$ of such a local fibration where $-\delta\leq t\leq \delta$ for some $0<\delta\leq \epsilon$. Then $N_t=\mathtt{exp}_{V_t}(N_0)$ are the desired $1$-parameter family of smooth special Lagrangian manifolds.
\end{proof}
\begin{remark} [on the proof] Strictly speaking, the domain of $F$ is not a Banach space because of the $(-\epsilon, \epsilon)$ part. This minor difficulty can be overcome by either using a cut-off function of $t$ or reparametrizing $t$ by a diffeomorphism between $(-\epsilon, \epsilon)$ and $\mathbf{R}$ preserving $0$.

\end{remark}


\begin{thebibliography}{1}

\bibitem{Besse}A. Besse, Einstein Manifolds, Springer-Verlag Berlin Heidelberg 1987 
\bibitem{Bryant1}R. Bryant, Metrics with exceptional holonomy, Ann. Math. 2, vol 126 (1987), 525-576 
\bibitem{Bryant2}R. Bryant, Minimal Lagrangian submanifolds of K\"{a}hler-Einstein Manifolds, Lecture Notes in Mathematics 1255, 1985
\bibitem{Bryant3}R. Bryant, Conformal and minimal immersions of compact surfaces into the 4-sphere, J. Differential Geometry 17 (1982), 455-473
\bibitem{Bryant4}R. Bryant, Second order families of special Lagrangian 3-folds, Centre de Recherches Math\'ematiques Precedings and Lecture Notes, vol 40 (2006), 63-98 
\bibitem{BGC} R. L. Bryant, S. S. Chern, R. B. Gardner, H. L. Goldschmidt, and P. A. Griffiths, Exterior Differential Systems, Springer-Verlag, New York, 1991
\bibitem{De-Ka} D. DeTurck and J. Kazdan, Some regularity theorems in Riemannian geometry, Ann. Scient. Ec. Norm. Sup. $4^{\circ}$ s\'erie, 14 (1981), 249-260 
\bibitem{Goldstein} E. Goldstein, Calibrated fibrations, math. DG/9911093, 1999
\bibitem{Harvey} R. Harvey, Spinors and Calibrations, Academic Press, Inc., 1990 
\bibitem{HarveyLawson} R. Harvey and H. B. Lawson, Calibrated geometries, Acta Mathematica 148 (1982), 47-157 
\bibitem{Joyce} D. Joyce, Lectures on Calabi-Yau and special Lagrangian geometry, math.DG/0108088, 2002
\bibitem{Joyce1} D. Joyce, On counting special Lagrangian homology 3-spheres, Contemp. Math. 314 (2002), 125-151
\bibitem{McLean} R. McLean, Deformations of calibrated submanifolds, Communications in Analysis and Geometry 6 (1998), 705-747
\bibitem{SYZ} A. Strominger, Shing-Tung Yau, and E. Zaslow, Mirror symmetry is T-duality, Nuclear Phys. B 479 (1996), no. 1-2, 243-259
\end{thebibliography}
\end{document}